\documentclass[11pt]{article}

\usepackage[a4paper,margin=1in]{geometry}
\usepackage[utf8]{inputenc}
\usepackage[T1]{fontenc}
\usepackage{lmodern}
\usepackage{amsmath,amssymb,amsthm,mathtools}
\usepackage{enumitem}
\usepackage{microtype}
\usepackage{graphicx}
\usepackage{tikz}
\usetikzlibrary{arrows.meta,decorations.markings}
\usepackage{float}
\usepackage{cite}
\usepackage{booktabs,tabularx}
\usepackage[hidelinks]{hyperref}

\newtheorem{theorem}{Theorem}[section]
\newtheorem{proposition}[theorem]{Proposition}
\newtheorem{lemma}[theorem]{Lemma}
\newtheorem{corollary}[theorem]{Corollary}
\newtheorem{remark}[theorem]{Remark}

\theoremstyle{definition}
\newtheorem{example}[theorem]{Example}
\numberwithin{equation}{section}

\newcommand{\R}{\mathbb{R}}
\newcommand{\ip}[2]{\langle #1,#2\rangle}
\newcommand{\norm}[1]{\|#1\|}
\newcommand{\conv}{\operatorname{conv}}
\newcommand{\inte}{\operatorname{int}}
\newcommand{\FR}{\mathrm{FR}}
\newcommand{\DY}{\mathrm{DY}}
\newcommand{\CD}{\mathrm{CD}}
\newcommand{\PRP}{\mathrm{PRP}}
\newcommand{\HS}{\mathrm{HS}}

\newcommand{\what}{\hat w}
\newcommand{\Kle}{\preceq_K}
\newcommand{\Klt}{\prec_K}
\newcommand{\D}{\mathcal{D}}
\newcommand{\doi}[1]{\url{https://doi.org/#1}}
\title{On the Fletcher--Reeves and Dai--Yuan Conjugate Gradient Methods for Vector
Optimization}
\author{
L. R. Lucambio P\'erez\thanks{Instituto de Matem\'atica e Estat\'istica,
Universidade Federal de Goi\'as, Avenida Esperan\c{c}a, s/n,
Campus Samambaia, Goi\^ania, GO 74690-900, Brazil
(\href{mailto:lrlp@ufg.br}{lrlp@ufg.br},
\href{mailto:lfprudente@ufg.br}{lfprudente@ufg.br}).}
\and
L. F. Prudente\footnotemark[1]
}
\date{}

\hypersetup{
  pdftitle={On the Fletcher--Reeves and Dai--Yuan Conjugate Gradient Methods for Vector Optimization},
  pdfauthor={L. R. Lucambio Perez and L. F. Prudente}
}

\begin{document}
\maketitle

\begin{abstract}
We study the descent and global convergence of the
Fletcher--Reeves (FR) and Dai--Yuan (DY) conjugate gradient
methods for vector optimization.
A biobjective example with strongly convex quadratic
objectives shows that the pure FR method may fail to
generate descent directions under the strong Wolfe
conditions, even with a curvature parameter smaller
than $1/2$.
For the pure DY method, the standard Wolfe conditions
yield descent but, unlike in the scalar case,
do not ensure global convergence.
We construct smooth biobjective examples showing that
both pure methods may fail to converge under the standard
Wolfe conditions, even when every search direction is a
descent direction and the initial level set is compact.
Nevertheless, suitable Wolfe-type line search conditions
guarantee sufficient descent and global convergence
of both pure methods.
For FR, we provide an upper bound on the curvature
tolerance, while for DY the strong Wolfe conditions suffice.
The key tool is a new estimate controlling the growth
of the search directions through the decrease
of the objective.
The convergence analysis covers scaled and hybrid parameters
and does not require bounded level sets.
We also establish the existence of Wolfe step sizes
for general ordering cones under a lower-bound
assumption on the objective values.


\end{abstract}

\section{Introduction}

Conjugate gradient (CG) methods constitute a class of unconstrained optimization algorithms which
combine simple recurrences with low storage requirements.
Their origins go back to the Hestenes--Stiefel (HS) method
\cite{HS52}, introduced in 1952 for the solution of linear
systems with symmetric positive definite matrices.
For nonlinear optimization, the first CG method was
the Fletcher--Reeves (FR) method \cite{FR64},
introduced in 1964.
Since then, numerous conjugate parameters have been proposed
and analyzed, including the Polak--Ribi\`ere--Polyak (PRP)
\cite{PR69,P69}, Dai--Yuan (DY) \cite{DY99}, and
Hager--Zhang (HZ) \cite{HZ05} formulas.
We refer to \cite{HZ06} for a comprehensive survey.

Unconstrained vector optimization generalizes scalar
minimization to a vector-valued objective
$F\colon\R^n\to\R^m$, whose values are compared through
a partial order induced by a cone $K\subset\R^m$.
In this paper, $F$ is continuously differentiable and
$K$ is closed, convex, and pointed, with nonempty interior.
We consider the problem
\begin{equation}\label{eq:prob}
\operatorname{minimize}_K\, F(x)
=\bigl(f_1(x),\ldots,f_m(x)\bigr)^{\top},
\qquad x\in\R^n.
\end{equation}
When $K=\R^m_+$, \eqref{eq:prob} is an unconstrained
multiobjective optimization problem.

Nonlinear CG methods for \eqref{eq:prob} were first
introduced in 2018 in \cite{LPP18}.
That work extended the FR, conjugate descent (CD), DY, PRP,
and HS parameters to vector optimization and established
vector counterparts of the Wolfe conditions and the
Zoutendijk condition. Since then, CG methods for vector and multiobjective
optimization have become an active research topic.
Subsequent developments include extensions of the HZ method
\cite{GP20,HZC24} and the Liu--Storey (LS) method \cite{GLP22},
as well as spectral, modified, and hybrid schemes
\cite{HCL23,YAI23,YK24,HLZC24,ZHCLL25,PZY26,CLLOZ26,HLC26}, and methods on
Riemannian manifolds \cite{NH23,CFS25}.

In the framework of \cite{LPP18}, the vector steepest
descent direction $v(x)$ of \cite{GDS05} replaces
the negative gradient.
Starting from $x^0\in\R^n$, these methods generate
\begin{equation}\label{eq:it}
x^{k+1}=x^k+\alpha_kd^k,
\end{equation}
where
\begin{equation}\label{eq:dir}
d^k=
\begin{cases}
v(x^k), & k=0,\\
v(x^k)+\beta_kd^{k-1}, & k\ge1.
\end{cases}
\end{equation}
Here, $\beta_k$ is the \emph{conjugate parameter} and $\alpha_k>0$
is chosen by a vector Wolfe line search.
The method stops if $v(x^k)=0$, which is equivalent
to $K$-criticality of $x^k$.
For an infinite sequence of iterates, global convergence
means
\[
\liminf_{k\to\infty}\norm{v(x^k)}=0.
\]
The definitions of $v$, $K$-criticality, the conjugate
parameters, and the line search conditions are recalled
in Section~\ref{sec:prel}.

We focus on the descent and global convergence of the \emph{pure}
FR and DY methods, corresponding to
$\beta_k=\beta_k^{\FR}$ and $\beta_k=\beta_k^{\DY}$,
respectively.
In scalar optimization, under standard assumptions,
the FR method generates sufficient descent directions
and is globally convergent with a strong Wolfe line search
whose curvature parameter satisfies $c_2<1/2$; see \cite{AlB85}.
Under the standard Wolfe conditions, the scalar FR
method is also globally convergent if every search direction
is a descent direction and the initial level set is bounded
\cite[Corollary~2]{Dai99}.
The DY method generates descent directions and is globally
convergent under the standard Wolfe conditions; see \cite{DY99}.
The distinction between descent and sufficient descent
in scalar CG convergence analysis was investigated
in \cite{DY96,DHLSYY00}.

For vector optimization, \cite{LPP18} established convergence
under suitable regularity and descent assumptions for
\[
0\le\beta_k\le\delta\beta_k^{\FR},\ \ 0\le\delta<1,
\qquad\text{or}\qquad
\beta_k=\eta\beta_k^{\DY},\ \ 0\le\eta<\frac{1-c_2}{1+c_2},
\]
where $c_2$ is the curvature parameter of the line search.
In \cite{LPP18}, the FR convergence result under the
standard Wolfe conditions assumes sufficient descent.
The result under the strong Wolfe conditions assumes
$K$-descent directions and a bounded level set.
The DY result in that work also uses the strong
Wolfe conditions.
The admissible range for the scaled DY parameter can
be enlarged to $0\le\eta<1$ as a consequence of
\cite[Theorem~3]{HLZC24}, but still excludes
the pure method.

The descent question for FR and the convergence questions
for the pure FR and DY methods were raised in \cite{LPP18}.
We organize them into three questions:
\begin{enumerate}[label=(Q\arabic*)]
\item Does the pure FR method with a strong Wolfe line search
and $c_2<1/2$ always generate $K$-descent directions?
\item Under what descent and line search conditions can global
convergence be established for $0\le\beta_k\le\beta_k^{\FR}$,
including the pure FR method?
\item Is the pure DY method globally convergent under
the strong Wolfe conditions, and do the standard Wolfe
conditions suffice?
\end{enumerate}

In this paper, we answer all three questions and provide
a unified analysis of the descent and global convergence
properties of the FR and DY methods in vector optimization.
Our convergence results also apply to several hybrid
choices of the conjugate parameter.

\begin{itemize}
\item The answer to (Q1) is negative.
We construct a biobjective problem with strongly convex
quadratic objectives for which the pure FR method loses
$K$-descent under a strong Wolfe line search
with $c_2<1/2$ (Example~\ref{ex:fr}).

\item For (Q2), we establish global convergence
for $0\le\beta_k\le\beta_k^{\FR}$ under sufficient descent
and the standard Wolfe conditions.
Under a generalized Wolfe line search that includes
the strong Wolfe conditions, $K$-descent alone suffices,
as in the scalar case \cite[Theorem~2.3]{DY96}.
We also give a line search condition that guarantees sufficient
descent and hence global convergence of the pure FR method
(Theorem~\ref{thm:fr}). Under the standard Wolfe conditions,
however, the pure method may fail to converge even when every
direction is a $K$-descent direction and the initial level set
is compact (Example~\ref{ex:frcompact}).

\item For (Q3), we show that the pure DY method generates
$K$-descent directions under the standard Wolfe conditions
and is globally convergent under the strong Wolfe conditions
(Theorem~\ref{thm:dy}). Under the standard Wolfe conditions,
global convergence may fail even on a smooth biobjective
problem with a compact initial level set (Example~\ref{ex:dy}).
\end{itemize}

We also establish the existence of Wolfe step sizes for
general ordering cones under the lower-bound assumption
used in our convergence analysis (Proposition~\ref{lem:exist}).
This removes the finite-generation restriction on the cone
used in \cite[Proposition~3.2]{LPP18}.

The main difficulty in the convergence analysis lies
in controlling the growth of the search directions.
In the scalar case, the Wolfe conditions control the inner
product $\ip{v(x^k)}{d^{k-1}}$ appearing in the expansion
of $\norm{d^k}^2$.
Their vector counterparts do not provide the same control.
To overcome this difficulty, we bound this inner product
from above in terms of the decrease of $F$ along
the preceding step (Lemma~\ref{lem:key}).
Together with the finite total decrease of $F$,
this estimate provides the key tool for establishing
a general convergence result for CG methods
(Theorem~\ref{thm:general}), which is then applied
to the FR and DY methods.
Our analysis assumes Lipschitz continuity of the Jacobian $JF$ of $F$ on the
convex hull of the initial level set and a lower-bound
condition for $K$-nonincreasing sequences of objective values.
Our convergence results do not require a bounded level set.

Section~\ref{sec:prel} presents the preliminaries and establishes
the existence of admissible Wolfe step sizes.
Section~\ref{sec:frcounterexample} shows that the pure FR method
may lose descent under the strong Wolfe conditions.
Section~\ref{sec:main} develops the general convergence framework,
which we apply to FR and DY in Sections~\ref{sec:fr}
and~\ref{sec:dy}, respectively.
Section~\ref{sec:nonconvergence} gives counterexamples to global
convergence of both methods under the standard Wolfe conditions
on problems with compact initial level sets.
Section~\ref{sec:final} presents our conclusions.


\section{Preliminaries and Auxiliary Results}\label{sec:prel}

This section collects the tools used throughout the paper.
Section~\ref{sec:basic} recalls the basic concepts of vector optimization and the
CG method of \cite{LPP18}, and Section~\ref{sec:wolfe} establishes the existence
of Wolfe step sizes under mild assumptions and recalls the Zoutendijk condition.

\subsection{Vector Optimization Concepts}\label{sec:basic}

The partial order induced by $K$ is given by $u\Kle v$ if $v-u\in K$, and the
corresponding strict order by $u\Klt v$ if $v-u\in\inte(K)$.
A point $x^*\in\R^n$ is \emph{$K$-optimal} if there is no $y\in\R^n$ such that
$F(y)\Kle F(x^*)$ and $F(y)\ne F(x^*)$, and \emph{weakly $K$-optimal} if there is no
$y\in\R^n$ such that $F(y)\Klt F(x^*)$. When $K=\R^m_+$, these are the notions of
\emph{Pareto} and \emph{weak Pareto optimality}.
We say that a direction $d\in\R^n$ is a \emph{$K$-descent direction} for $F$ at $x$ if
$JF(x)d\Klt 0$, and $x$ is \emph{$K$-critical} if there is no such direction.
$K$-criticality is a necessary condition for both notions of optimality \cite{GDS05}.

Let $K^*=\{w\in\R^m\mid\ip{w}{y}\ge0\ \forall y\in K\}$ be the positive polar cone
of $K$, and let $C\subset K^*$ be a compact set such that $\norm{w}=1$ for all
$w\in C$ and 
\[K^*=\operatorname{cone}(\conv(C)),\]
i.e., $K^*$ is the conic hull of the convex hull of $C$.
In multiobjective optimization,
$K=K^*=\R^m_+$ and $C$ can be taken as the canonical basis of $\R^m$. More
generally, if $K$ is polyhedral, then so is $K^*$, and $C$ can be taken as the
finite set of unit vectors spanning the extreme rays of $K^*$.

Define $\D\colon\R^n\times\R^n\to\R$ by
\[
\D(x,d)=\max_{w\in C}\ip{JF(x)d}{w},
\]
where $JF(x)$ is the Jacobian of $F$ at $x$. The function $\D$ characterizes $K$-descent directions and $K$-critical
points: $d$ is a $K$-descent direction for $F$ at $x$ if and only if $\D(x,d)<0$, and
$x$ is $K$-critical if and only if $\D(x,d)\ge0$ for all $d\in\R^n$.
The next lemma collects two further properties of $\D$, which follow directly from
its definition and from $\norm{w}=1$ for $w\in C$.

\begin{lemma}\label{lem:D}
\begin{enumerate}[label=\textup{(\alph*)}]
\item For each $x\in\R^n$, $\D(x,\cdot)$ is sublinear, that is,
 $\D(x,td)=t\D(x,d)$ and $\D(x,d+\bar d)\le\D(x,d)+\D(x,\bar d)$ for all
 $d,\bar d\in\R^n$ and $t\ge0$.
\item If $JF$ is Lipschitz continuous with constant $L$ on a set $S\subset\R^n$, then
 $|\D(x,d)-\D(y,d)|\le L\norm{x-y}\norm{d}$ for all $x,y\in S$ and $d\in\R^n$.
\end{enumerate}
\end{lemma}

The next function extends the steepest descent direction to the vector setting
\cite{GDS05}. Define $v\colon\R^n\to\R^n$ by
\[
v(x)=\operatorname*{arg\,min}_{d\in\R^n}\Bigl\{\D(x,d)+\tfrac12\norm{d}^2\Bigr\}.
\]
Since $\D(x,\cdot)$ is convex, the objective is strongly convex, so $v(x)$ exists and
is unique.

\begin{lemma}[{\cite[Lemma~3.3 and Section~5]{GDS05}}]\label{lem:basic}
Let $x\in\R^n$.
\begin{enumerate}[label=\textup{(\alph*)}]
\item $\D(x,v(x))=-\norm{v(x)}^2$. In particular, $x$ is $K$-critical if and only if
 $v(x)=0$, and otherwise $v(x)$ is a $K$-descent direction for $F$ at $x$.
\item There exists $w(x)\in\conv(C)$ such that $v(x)=-JF(x)^{\top}w(x)$.
\item The mapping $v$ is continuous.
\end{enumerate}
\end{lemma}

By Lemma~\ref{lem:basic}(a), $\D(x,v(x))=-\norm{v(x)}^2$. Motivated by this identity,
we say that the search directions $d^k$ in \eqref{eq:it} satisfy the sufficient
descent condition \cite{LPP18} if there exists $\kappa>0$ such that
\begin{equation}\label{eq:sd}
\D(x^k,d^k)\le-\kappa\norm{v(x^k)}^2\qquad\text{for all }k\ge0.
\end{equation}
In turn, Lemma~\ref{lem:basic}(b) yields an inequality that plays a central role in
our analysis. Since a linear function has the same maximum over $C$ and over
$\conv(C)$, we have, for every $d\in\R^n$,
\begin{equation}\label{eq:vd}
\ip{v(x)}{d}=-\ip{w(x)}{JF(x)d}\ge-\D(x,d).
\end{equation}
In the scalar case, \eqref{eq:vd} holds with equality. In the vector case, the
inequality may be strict. In multiobjective optimization, for instance, $\D(x,d)$ is
the largest directional derivative, whereas
$\ip{v(x)}{d}$ is minus a weighted average of the directional derivatives of the
objectives that are active at $v(x)$. This gap is the source of the difficulties
addressed in this paper.

We now turn to the CG method \eqref{eq:it}--\eqref{eq:dir}. It stops if
$v(x^k)=0$, since then $x^k$ is $K$-critical. From now on, we assume that this does
not happen, that is, $v(x^k)\ne0$ for all $k$. Using
Lemma~\ref{lem:basic}(a), the vector extensions
of the FR and DY parameters of \cite{LPP18} read
\[
\beta_k^{\FR}=\frac{\norm{v(x^k)}^2}{\norm{v(x^{k-1})}^2},\qquad
\beta_k^{\DY}=\frac{\norm{v(x^k)}^2}{\D(x^k,d^{k-1})-\D(x^{k-1},d^{k-1})}.
\]
Moreover, if $\beta_k\ge0$, applying Lemma~\ref{lem:D}(a) to \eqref{eq:dir} gives the basic inequality
\begin{equation}\label{eq:basic}
\D(x^k,d^k)\le-\norm{v(x^k)}^2+\beta_k\D(x^k,d^{k-1}),
\end{equation}
which holds with equality in the scalar case.

\subsection{Wolfe Step Sizes and the Zoutendijk Condition}\label{sec:wolfe}

In this subsection, $\{x^k\}$ is generated by \eqref{eq:it}, where $d^k$ is any
$K$-descent direction for $F$ at $x^k$, not necessarily given by \eqref{eq:dir},
and $\alpha_k>0$ is a step size satisfying Wolfe-type conditions.
Let $e\in\inte(K)$ be such that
\begin{equation}\label{eq:e}
0<\ip{w}{e}\le1\qquad\text{for all }w\in C.
\end{equation}
Such a vector exists, see \cite[Definition~3.1]{LPP18}.
Given $0<c_1<c_2<1$, we say that $\alpha_k>0$ satisfies the Wolfe conditions
with tolerance $\varepsilon_k\in[0,+\infty]$ if
\begin{gather}
F(x^k+\alpha_kd^k)\Kle F(x^k)+c_1\,\alpha_k\D(x^k,d^k)\,e,\label{eq:w1}\\
c_2 \D(x^k,d^k)\le\D(x^k+\alpha_kd^k,d^k)\le\varepsilon_k.\label{eq:w2}
\end{gather}
Condition \eqref{eq:w1} is known as the
\emph{Armijo (or sufficient decrease) condition},
whereas \eqref{eq:w2} is a
\emph{generalized curvature condition}.
The choice $\varepsilon_k=+\infty$ gives the \emph{standard Wolfe conditions} and $\varepsilon_k=-c_2\D(x^k,d^k)$ gives the \emph{strong Wolfe conditions}; see \cite{LPP18}.

Let $\mathcal L=\{x\in\R^n\mid F(x)\Kle F(x^0)\}$. Our analysis relies on the
following assumptions.
\begin{enumerate}[label=\textup{(A\arabic*)}]
\item $JF$ is Lipschitz continuous with constant $L>0$ on $\conv(\mathcal L)$.
\item Every sequence in $F(\mathcal L)$ that is nonincreasing with respect to
 $\Kle$ is bounded below with respect to $\Kle$.
\end{enumerate}

Assumption (A1) holds, for instance, when $\mathcal L$ is bounded and $JF$ is
locally Lipschitz continuous. In \cite{LPP18}, $JF$ is required to be Lipschitz
continuous on an open set containing $\mathcal L$. The proofs below use (A1) only
along the segments $[x^k,x^{k+1}]$ joining successive iterates, which are contained
in $\conv(\mathcal L)$. Assumption (A2) holds,
for instance, if $F$ is bounded below on $\mathcal L$, that is, if there is
$\bar F\in\R^m$ with $\bar F\Kle F(x)$ for all $x\in\mathcal L$. In multiobjective
optimization, this sufficient condition means that each $f_i$ is bounded below on
$\mathcal L$.

In \cite[Proposition~3.2]{LPP18}, the existence of step sizes satisfying the Wolfe
conditions was established for finitely generated cones $K$, under a lower bound on
$F$ along the search ray. A similar lower bound was used in
\cite[Theorem~3.1]{GRPY25} for set optimization with a general cone. The next proposition
shows that, under (A2), such step sizes exist without requiring $K$ to be finitely
generated.

\begin{proposition}\label{lem:exist}
Assume that \textup{(A2)} holds, let $x\in\mathcal L$, and let $d$ be a $K$-descent
direction for $F$ at $x$. Then there is a nonempty open interval
$I\subset(0,\infty)$ such that every $\alpha\in I$ satisfies
\[
F(x+\alpha d)\Kle F(x)+c_1\alpha\D(x,d)e
\qquad\text{and}\qquad
c_2\D(x,d)<\D(x+\alpha d,d)<0 .
\]
In particular, every $\alpha\in I$ satisfies \eqref{eq:w1}--\eqref{eq:w2} at $x$
along $d$ for every tolerance $\varepsilon\in[0,+\infty]$.
\end{proposition}
\begin{proof}
Fix $c\in(c_1,c_2)$, and define $h\colon[0,+\infty)\to\R$ by
\[
h(u):=\D(x+ud,d).
\]
Then $h$ is continuous and $h(0)=\D(x,d)<c\D(x,d)$, since $\D(x,d)<0$. Recall
also that, since $K^*=\operatorname{cone}(\conv(C))$ and $K$ is closed and convex,
$y\in K$ if and only if $\ip{w}{y}\ge0$ for all $w\in C$.

We first show that $F$ decreases along every interval on which
$h$ is bounded above by $c\D(x,d)$.
Fix $0\le s<t$ and suppose that
\[
h(u)\le c\D(x,d)\qquad\text{for every }u\in[s,t].
\]
For any $w\in C$, we have
$\ip{w}{JF(x+ud)d}\le h(u)$ for every $u\in[s,t]$.
Hence,
\[
\begin{aligned}
\ip{w}{F(x+td)-F(x+sd)}
&=\int_s^t\ip{w}{JF(x+ud)d}\,du\\
&\le c(t-s)\D(x,d)
\le c(t-s)\D(x,d)\ip{w}{e},
\end{aligned}
\]
where the last inequality follows from $\D(x,d)<0$ and \eqref{eq:e}. Since
$w\in C$ is arbitrary,
\begin{equation}\label{eq:decr}
F(x+td)\Kle F(x+sd)+c(t-s)\D(x,d)e .
\end{equation}
In particular, since $c\D(x,d)e\in-K$, if $h\le c\D(x,d)$ on $[0,t]$, then
\eqref{eq:decr} with $s=0$ gives $F(x+td)\Kle F(x)\Kle F(x^0)$, so
$x+td\in\mathcal L$.

Next, we claim that $h$ cannot stay below $c\D(x,d)$ forever. Indeed, suppose,
to the contrary, that $h(t)\le c\D(x,d)$ for all $t\ge0$. For every integer
$j\ge0$, \eqref{eq:decr} with $s=j$ and $t=j+1$ gives $F(x+(j+1)d)\Kle F(x+jd)$.
Hence $\{F(x+jd)\}_{j\ge0}\subset F(\mathcal L)$ is nonincreasing with respect to
$\Kle$. By (A2), there is $\bar F$ with $\bar F\Kle F(x+jd)$ for all $j$. Fixing
$w_0\in C$ and using \eqref{eq:decr} with $s=0$ and $t=j$, we obtain
\[
\ip{w_0}{\bar F}\le\ip{w_0}{F(x+jd)}\le\ip{w_0}{F(x)}+j c\D(x,d)\ip{w_0}{e}
\to-\infty,
\]
which is a contradiction. Thus the set $S:=\{t>0\mid h(t)>c\D(x,d)\}$ is
nonempty. Setting $t_*:=\inf S$, continuity yields $t_*>0$, $h(t_*)=c\D(x,d)$,
and $h(t)\le c\D(x,d)$ for all $t\in[0,t_*]$.

We now show that both inequalities in the proposition hold strictly at $t_*$. By
\eqref{eq:decr} with $s=0$ and $t=t_*$,
\[
F(x)+c t_*\D(x,d)e-F(x+t_*d)\in K.
\]
Since $c_1<c$, $\D(x,d)<0$, $t_*>0$, and $e\in\inte(K)$, we have
\[
(c_1-c)t_*\D(x,d)e\in\inte(K).
\]
Consequently, adding these two terms and using $K+\inte(K)\subset\inte(K)$, we
obtain
\[
F(x)+c_1 t_*\D(x,d)e-F(x+t_*d)\in\inte(K).
\]
Moreover, since $c<c_2$ and $\D(x,d)<0$,
\[
c_2\D(x,d)<h(t_*)=c\D(x,d)<0.
\]
Since $F$ and $h$ are continuous and $\inte(K)$ is open,
both strict conditions persist on a sufficiently small
open interval $I\subset(0,\infty)$ centered at $t_*$. Every
$\alpha\in I$ therefore satisfies both inequalities in the proposition, which completes
the proof.
\end{proof}

Proposition~\ref{lem:exist} guarantees the existence of
suitable step sizes but does not provide a procedure
for computing them. For finitely generated cones,
the line search algorithm of \cite{LPP19} provides
a natural basis for incorporating the upper curvature
tolerances considered here.

We close this
section by recalling that the Zoutendijk condition holds for any descent method
with Wolfe step sizes. Its proof in \cite{LPP18} remains valid under (A1)--(A2).

\begin{proposition}[{\cite[Proposition~3.3]{LPP18}}]
\label{lem:zout}
Assume that \textup{(A1)} and \textup{(A2)} hold,
and consider \eqref{eq:it}, where $d^k$ is a
$K$-descent direction for $F$ at $x^k$ and
$\alpha_k$ satisfies \eqref{eq:w1}--\eqref{eq:w2}
for all $k$, with arbitrary tolerances
$\varepsilon_k\in[0,+\infty]$. Then
\[
\sum_{k\ge0}
\frac{\D(x^k,d^k)^2}{\norm{d^k}^2}<\infty.
\]
\end{proposition}


\section{The Fletcher--Reeves Method May Fail to Generate Descent Directions}\label{sec:frcounterexample}

Throughout this section, we consider the method \eqref{eq:it}--\eqref{eq:dir} with
$0\le\beta_k\le\beta_k^{\FR}$, and we assume that $\alpha_k$ satisfies the strong
Wolfe conditions. Let 
\[r_k:=-\frac{\D(x^k,d^k)}{\norm{v(x^k)}^2},\] 
so that $d^k$ is a
$K$-descent direction if and only if $r_k>0$, and \eqref{eq:sd} means $r_k\ge\kappa$
for all $k$. 

If $d^{k-1}$ is a $K$-descent direction, then \eqref{eq:basic} and
\eqref{eq:w2} with $\varepsilon_{k-1}=-c_2\D(x^{k-1},d^{k-1})$ give
\[
\D(x^k,d^k)\le-\norm{v(x^k)}^2+\beta_k^{\FR}c_2|\D(x^{k-1},d^{k-1})|
=-\norm{v(x^k)}^2(1-c_2 r_{k-1}),
\]
that is,
\begin{equation}\label{eq:rk}
r_k\ge1-c_2 r_{k-1}.
\end{equation}
Since $d^0=v(x^0)$, Lemma~\ref{lem:basic}(a) gives $r_0=1$, and hence
$\D(x^1,d^1)\le-(1-c_2)\norm{v(x^1)}^2$. Thus $d^1$ is always a $K$-descent
direction, with $r_1\ge1-c_2$.

In the scalar case, \eqref{eq:basic} is an equality, and the same argument also
gives $r_k\le1+c_2 r_{k-1}$. Since $r_0=1$, an induction gives
$r_k\le1/(1-c_2)$ and $r_k\ge(1-2c_2)/(1-c_2)>0$ for all $k$ when
$c_2<1/2$. This is Al-Baali's argument \cite{AlB85}, which shows that the
scalar FR method with a strong Wolfe line search and $c_2<1/2$ satisfies the
sufficient descent condition with $\kappa=(1-2c_2)/(1-c_2)$.

In the vector case, \eqref{eq:basic} is only an inequality, and no upper bound on
$r_{k-1}$ is available. As soon as $r_{k-1}>1/c_2$, the lower bound \eqref{eq:rk} becomes
negative and no longer guarantees that $d^k$ is a descent direction. The following
example shows that descent can indeed be lost, already for $k=2$.

\begin{example}\label{ex:fr}
Let $K=\R^2_+$, $C$ the canonical basis of $\R^2$, $e=(1,1)^{\top}$, $c_1=10^{-4}$,
$c_2=0.1$, $x^0=(0,0)^{\top}$, and $F=(f_1,f_2)$ with
\[
f_1(x)=\tfrac12(x_1-1)^2+1000\,x_1x_2+10^6x_2^2,\qquad
f_2(x)=\tfrac12(x_1-21)^2+\tfrac12(x_2+10)^2 .
\]
Both functions are strongly convex quadratics, so \textup{(A1)} and \textup{(A2)}
hold. We have 
\[\nabla f_1(x)=(x_1-1+1000x_2,\ 1000x_1+2\cdot10^6x_2)^{\top} \quad \text{and}\quad \nabla f_2(x)=(x_1-21,\ x_2+10)^{\top}.\] 
For $m=2$, $-v(x)$ is the minimal-norm element of the segment
$[\nabla f_1(x),\nabla f_2(x)]$ \cite{FS00}. An endpoint $a$ of a segment $[a,b]$ is
its minimal-norm element if and only if $\ip{a}{b-a}\ge0$. Hence
\[
v(x)=\begin{cases}
-\nabla f_1(x) & \text{if } \ip{\nabla f_1(x)}{\nabla f_2(x)-\nabla f_1(x)}\ge0,\\
-\nabla f_2(x) & \text{if } \ip{\nabla f_2(x)}{\nabla f_1(x)-\nabla f_2(x)}\ge0.
\end{cases}
\]
Otherwise, the minimal-norm element lies in the interior of the segment, both
objectives are active, and $\ip{\nabla f_i(x)}{v(x)}=-\norm{v(x)}^2$ for $i=1,2$.

\emph{Iteration $0$.} $\nabla f_1(x^0)=(-1,0)^{\top}$ and
$\nabla f_2(x^0)=(-21,10)^{\top}$. Since
\[
\ip{\nabla f_1(x^0)}{\nabla f_2(x^0)-\nabla f_1(x^0)}=20\ge0,
\]
we have
$d^0=v(x^0)=-\nabla f_1(x^0)=(1,0)^{\top}$ and $\D(x^0,d^0)=-1$. The step size
$\alpha_0=1$ gives 
\[x^1=(1,0)^{\top}, \quad \nabla f_1(x^1)=(0,1000)^{\top}, \quad \nabla f_2(x^1)=(-20,10)^{\top},\] 
and $\D(x^1,d^0)=\max\{0,-20\}=0$. Moreover,
$f_1(x^1)-f_1(x^0)=-0.5$ and $f_2(x^1)-f_2(x^0)=-20.5$. Hence $\alpha_0$ satisfies
the strong Wolfe conditions (it is an exact step).

\emph{Iteration $1$.} Since
\[
\ip{\nabla f_2(x^1)}{\nabla f_1(x^1)-\nabla f_2(x^1)}=9500\ge0,
\]
we have
$v(x^1)=-\nabla f_2(x^1)=(20,-10)^{\top}$, $\beta_1^{\FR}=500$, $d^1=(520,-10)^{\top}$, and
\[
\D(x^1,d^1)=\max\{-10000,\,-10500\}=-10000=-20\norm{v(x^1)}^2 ,
\]
that is, $r_1=20>1/c_2$. Take $\alpha_1=5.6\times10^{-5}$. Then
$x^2=x^1+\alpha_1d^1=(1.02912,\,-0.00056)^{\top}$,
\[
\nabla f_1(x^2)=(-0.53088,\ -90.88)^{\top},\qquad
\nabla f_2(x^2)=(-19.97088,\ 9.99944)^{\top},
\]
so that
\[
\ip{\nabla f_1(x^2)}{d^1}=632.7424,\qquad\ip{\nabla f_2(x^2)}{d^1}=-10484.852.
\]
Hence $\D(x^2,d^1)=632.7424\le c_2|\D(x^1,d^1)|=1000$. Since $f_i$ is quadratic,
$s\mapsto\ip{\nabla f_i(x^1+sd^1)}{d^1}$ is affine, so the trapezoidal rule is
exact:
\[
f_i(x^2)-f_i(x^1)=\frac{\alpha_1}{2}\big(\ip{\nabla f_i(x^1)}{d^1}
+\ip{\nabla f_i(x^2)}{d^1}\big).
\]
This gives
\[
f_1(x^2)-f_1(x^1)=-4683.6288\,\alpha_1,\qquad
f_2(x^2)-f_2(x^1)=-10492.426\,\alpha_1,
\]
both smaller than $c_1\alpha_1\D(x^1,d^1)=-\alpha_1$. Hence $\alpha_1$ satisfies
the strong Wolfe conditions.

\emph{Iteration $2$.} We have
\[
\ip{\nabla f_1(x^2)}{\nabla f_2(x^2)}=0.53088\cdot19.97088-90.88\cdot9.99944
<11-900<0,
\]
so neither of the two cases in the formula for $v(x)$ above applies. Thus both objectives are active at $x^2$. Moreover, $v(x^2)\ne0$, because
every point of the segment $[\nabla f_1(x^2),\nabla f_2(x^2)]$ has a negative
first component. With
$\beta_2=\beta_2^{\FR}=\norm{v(x^2)}^2/\norm{v(x^1)}^2=\norm{v(x^2)}^2/500$, we
obtain
\[
\begin{aligned}
\D(x^2,d^2)&\ge\ip{\nabla f_1(x^2)}{v(x^2)+\beta_2d^1}
=-\norm{v(x^2)}^2+\frac{\norm{v(x^2)}^2}{500}\cdot632.7424\\
&=0.2654848\,\norm{v(x^2)}^2\approx83.34>0 .
\end{aligned}
\]
Therefore $d^2$ is not a $K$-descent direction for $F$ at $x^2$. Since
$\ip{\nabla f_2(x^2)}{d^2}=-\norm{v(x^2)}^2+\beta_2\ip{\nabla f_2(x^2)}{d^1}<0$,
equality holds above.

The failure is not tied to this particular step size. Indeed, for $\alpha\ge0$,
\[
\D(x^1+\alpha d^1,d^1)
=\ip{\nabla f_1(x^1+\alpha d^1)}{d^1}
=-10000+189870400\,\alpha,
\]
so the strong Wolfe curvature condition at iteration~$1$ holds precisely for
\[
\alpha_1\in\left[\frac{9000}{189870400},\frac{11000}{189870400}\right],
\]
and the Armijo condition holds throughout this interval. Moreover, for every
$\alpha_1$ in its upper quarter $[10500,11000]/189870400$,
both objectives are active at $x^2$ (their gradients have
negative inner product) and
$\ip{\nabla f_1(x^2)}{d^1}\ge500=\norm{v(x^1)}^2$.
Consequently,
\[
\D(x^2,d^2)\ge\ip{\nabla f_1(x^2)}{v(x^2)+\beta_2d^1}\ge-\norm{v(x^2)}^2
+\frac{\norm{v(x^2)}^2}{500}\cdot500=0.
\]
Hence $d^2$ is not a $K$-descent direction for any such $\alpha_1$.
\hfill$\diamond$
\end{example}

Example~\ref{ex:fr} answers (Q1): in the vector setting,
the strong Wolfe conditions with $c_2<1/2$
do not guarantee that the FR method generates descent
directions, even for strongly convex quadratic objectives.
Nevertheless, sufficient descent can be guaranteed by
imposing an appropriate upper curvature tolerance,
as shown in Theorem~\ref{thm:fr}(a).

\section{A Bound for the Cross Term and a General Convergence Result}\label{sec:main}

In this section, $\{x^k\}$ is generated by an iteration of the form
\eqref{eq:it}, where each $d^k$ is a $K$-descent direction for $F$ at $x^k$ and
each $\alpha_k$ satisfies \eqref{eq:w1}--\eqref{eq:w2}, with arbitrary tolerances
$\varepsilon_k\in[0,+\infty]$, and $v(x^k)\ne0$ for all $k$. We first collect some
estimates involving the decrease of $F$ along the iterates.

\begin{lemma}\label{lem:step}
Assume that \textup{(A1)} holds, and consider \eqref{eq:it} with $K$-descent
directions $d^k$ and step sizes satisfying \eqref{eq:w1}--\eqref{eq:w2} with
arbitrary tolerances $\varepsilon_k\in[0,+\infty]$. Let
$\theta:=\max_{w\in C}\ip{w}{e}\in(0,1]$ and
\begin{equation}\label{eq:delta}
\delta_k:=\norm{F(x^{k-1})-F(x^k)},\qquad k\ge1.
\end{equation}
Then:
\begin{enumerate}[label=\textup{(\alph*)}]
\item $\delta_k\ge c_1\theta\,\alpha_{k-1}|\D(x^{k-1},d^{k-1})|$ for all
 $k\ge1$, and $x^k\in\mathcal L$ for all $k\ge0$;
\item $(1-c_2)|\D(x^k,d^k)|\le L\,\alpha_k\norm{d^k}^2$ for all $k\ge0$;
\item if \textup{(A2)} also holds, then $\sum_{k\ge1}\delta_k<\infty$.
\end{enumerate}
\end{lemma}
\begin{proof}
(a) By \eqref{eq:w1},
$F(x^{k-1})-F(x^k)-c_1\alpha_{k-1}|\D(x^{k-1},d^{k-1})|e\in K$. Taking the inner
product with $w\in C\subset K^*$ such that $\ip{w}{e}=\theta$, and using
$\norm{w}=1$, gives 
\[\delta_k\ge\ip{w}{F(x^{k-1})-F(x^k)}\ge c_1\theta\alpha_{k-1}|\D(x^{k-1},d^{k-1})|.\]
Moreover, $F(x^k)\Kle F(x^{k-1})$
for all $k\ge1$, so $x^k\in\mathcal L$ for all $k$.

(b) By \eqref{eq:w2},
$\D(x^{k+1},d^k)-\D(x^k,d^k)\ge(1-c_2)|\D(x^k,d^k)|$, whereas (a), (A1), and
Lemma~\ref{lem:D}(b) give
$\D(x^{k+1},d^k)-\D(x^k,d^k)\le L\alpha_k\norm{d^k}^2$.
Combining these two inequalities proves (b).

(c) Since $K$ is closed, convex, and pointed, $\inte(K^*)\ne\emptyset$
\cite[Corollary~14.6.1]{Roc70}. Take $w_0\in\inte(K^*)$ and $r>0$ such that
$w_0-ru\in K^*$ whenever $\norm{u}\le1$, and set $\what:=w_0/r\in K^*$. For
$y\in K\setminus\{0\}$, taking $u=y/\norm{y}$ gives $\ip{w_0-ry/\norm{y}}{y}\ge0$,
and hence
\[
\norm{y}\le\ip{\what}{y}\qquad\text{for all }y\in K
\]
(for $K=\R^m_+$ one may take $\what=(1,\dots,1)^{\top}$). 
As shown in (a), $F(x^k)\Kle F(x^{k-1})$ for all $k\ge1$
and $x^k\in\mathcal L$ for all $k\ge0$.
Thus, by (A2), there exists $\bar F\in\R^m$ such that
$\bar F\Kle F(x^k)$ for all $k$. Since
$F(x^{k-1})-F(x^k)\in K$, applying the norm bound above with
$y=F(x^{k-1})-F(x^k)$ gives
$\delta_k\le\ip{\what}{F(x^{k-1})-F(x^k)}$.
Therefore, since $\what\in K^*$, for every $N\ge1$,
\[
\sum_{k=1}^N\delta_k
\le \sum_{k=1}^N
\ip{\what}{F(x^{k-1})-F(x^k)}
=\ip{\what}{F(x^0)-F(x^N)}
\le\ip{\what}{F(x^0)-\bar F}.
\]
Since $\delta_k\ge0$ and the right-hand side is independent
of $N$, the series converges, proving (c).
\end{proof}

Our convergence analysis relies on controlling the growth of
the search directions. For $k\ge1$, we have
\[
\norm{d^k}^2
=\norm{v(x^k)}^2
+2\beta_k\ip{v(x^k)}{d^{k-1}}
+\beta_k^2\norm{d^{k-1}}^2.
\]
In the scalar case,
$\ip{v(x^k)}{d^{k-1}}=-\D(x^k,d^{k-1})$,
so the curvature condition provides an upper bound for the
cross term. In the vector case, \eqref{eq:vd} gives only
the lower bound
$\ip{v(x^k)}{d^{k-1}}\ge-\D(x^k,d^{k-1})$.
Applying the same inequality to $-d^{k-1}$ gives
\[
\ip{v(x^k)}{d^{k-1}}\le\D(x^k,-d^{k-1})=-\min_{w\in C}\ip{JF(x^k)d^{k-1}}{w},
\]
so it suffices to bound $\D(x^k,-d^{k-1})$ from above.
The next lemma bounds $\D(z,-d^{k-1})$ from above uniformly for
$z\in[x^{k-1},x^k]$, in terms of the decrease of $F$ over the step.
Equivalently, it provides a common lower bound
for the directional derivatives along $d^{k-1}$ of the
scalarizations $\ip{w}{F(\cdot)}$, $w\in C$. In particular, it
bounds both $\ip{v(x^k)}{d^{k-1}}$ and $\ip{v(x^{k-1})}{d^{k-1}}$ from above.

\begin{lemma}\label{lem:key}
Assume that \textup{(A1)} holds, that $x^{k-1}\in\mathcal L$, that $d^{k-1}$ is a
$K$-descent direction for $F$ at $x^{k-1}$, and that $\alpha_{k-1}$ satisfies the
Wolfe conditions \eqref{eq:w1}--\eqref{eq:w2} with an arbitrary tolerance
$\varepsilon_{k-1}\in[0,+\infty]$. Then, for every $z\in[x^{k-1},x^k]$,
\[
\D(z,-d^{k-1})\ \le\ M\,\delta_k\,\frac{\norm{d^{k-1}}^2}{|\D(x^{k-1},d^{k-1})|},
\qquad
M:=L\Bigl[\frac{1}{1-c_2}+\frac{1}{2c_1\theta}\Bigr],
\]
where $\delta_k$ is given by \eqref{eq:delta}. In particular,
\[
\max\bigl\{
\ip{v(x^k)}{d^{k-1}},
\ip{v(x^{k-1})}{d^{k-1}}
\bigr\}
\le M\,\delta_k\,
\frac{\norm{d^{k-1}}^2}
     {|\D(x^{k-1},d^{k-1})|}.
\]
\end{lemma}
\begin{proof}
By \eqref{eq:w1}, the segment $[x^{k-1},x^k]$ is contained
in $\conv(\mathcal L)$. Let $z\in[x^{k-1},x^k]$.
By the fundamental theorem of calculus and (A1),
\[
\begin{aligned}
&\Bigl\|F(x^k)-F(x^{k-1})
-\alpha_{k-1}JF(z)d^{k-1}\Bigr\|
\\
&\quad=\Bigl\|\int_0^{\alpha_{k-1}}
\bigl[JF(x^{k-1}+td^{k-1})-JF(z)\bigr]
d^{k-1}\,dt\Bigr\|\\
&\quad\le
L\norm{d^{k-1}}
\int_0^{\alpha_{k-1}}
\norm{x^{k-1}+td^{k-1}-z}\,dt\\
&\quad\le\frac L2\alpha_{k-1}^2\norm{d^{k-1}}^2,
\end{aligned}
\]
where the last inequality follows because the integral is
convex as a function of $z$ and equals
$\alpha_{k-1}^2\norm{d^{k-1}}/2$ at both endpoints
of the segment.
Hence, for every $w\in C$,
using $\norm{w}=1$ and \eqref{eq:delta},
\[
\begin{aligned}
-\alpha_{k-1}\langle JF(z)d^{k-1},w\rangle
&=\langle F(x^{k-1})-F(x^k),w\rangle\\
&\quad+\langle F(x^k)-F(x^{k-1})
-\alpha_{k-1}JF(z)d^{k-1},w\rangle\\
&\le \|F(x^{k-1})-F(x^k)\|
+\frac L2\alpha_{k-1}^2\|d^{k-1}\|^2\\
&=\delta_k+\frac L2\alpha_{k-1}^2\|d^{k-1}\|^2.
\end{aligned}
\]
Dividing by $\alpha_{k-1}>0$ and taking the maximum over
$w\in C$ gives
\[
\D(z,-d^{k-1})
\le\frac{\delta_k}{\alpha_{k-1}}
+\frac L2\,\alpha_{k-1}\norm{d^{k-1}}^2.
\]
To estimate the two terms on the right-hand side,
Lemma~\ref{lem:step}(b) and (a) give
\[
\frac1{\alpha_{k-1}}
\le
\frac{L\norm{d^{k-1}}^2}
     {(1-c_2)|\D(x^{k-1},d^{k-1})|}
\qquad\text{and}\qquad
\alpha_{k-1}
\le
\frac{\delta_k}
     {c_1\theta|\D(x^{k-1},d^{k-1})|}.
\]
Substituting these two bounds proves the first claim.
Finally, \eqref{eq:vd} applied to the direction $-d^{k-1}$
gives
\[
\ip{v(z)}{d^{k-1}}\le\D(z,-d^{k-1}),
\]
which proves the last claim by taking
$z=x^{k-1}$ and $z=x^k$.
\end{proof}

The next theorem uses Lemma~\ref{lem:key} to establish
a general convergence result for CG methods. Its proof
exploits the finite total decrease of $F$, established
under (A2) in Lemma~\ref{lem:step}(c), to control the growth
of the search directions. The hypotheses are formulated
in terms of an auxiliary positive sequence $\{\omega_k\}$.
Suitable choices of this sequence will yield the FR and DY
convergence results in Sections~\ref{sec:fr} and~\ref{sec:dy}.

\begin{theorem}\label{thm:general}
Assume that \textup{(A1)} and \textup{(A2)} hold. Consider the method \eqref{eq:it}--\eqref{eq:dir} with
$\beta_k\ge0$, where $d^k$ is a $K$-descent direction for $F$ at $x^k$ and $\alpha_k$
satisfies the Wolfe conditions \eqref{eq:w1}--\eqref{eq:w2}, with arbitrary
tolerances $\varepsilon_k\in[0,+\infty]$, for all $k$.
Suppose that there exist constants $a,b>0$ and positive numbers $\omega_k$ such
that, for all $k$,
\[
\textup{(i)}\ \beta_k\le\frac{\omega_k}{\omega_{k-1}}\ (k\ge1),\quad
\textup{(ii)}\ \omega_k\ge a\norm{v(x^k)}^2,\quad
\textup{(iii)}\ -\D(x^k,d^k)\ge b\,\omega_k .
\]
Then $\liminf_{k\to\infty}\norm{v(x^k)}=0$.
\end{theorem}
\begin{proof}
Assume, by contradiction, that $\norm{v(x^k)}\ge\gamma>0$ for all $k$, and let
$s_k:=\norm{d^k}^2/\omega_k^2$. By \eqref{eq:w1}, $\{F(x^k)\}$ is $K$-nonincreasing,
so $x^k\in\mathcal L$ for all $k$ and Lemma~\ref{lem:key} applies at every
$k\ge1$. Let $k\ge1$. Since $\beta_k\ge0$, Lemma~\ref{lem:key} gives
\[
\norm{d^k}^2\le\norm{v(x^k)}^2
+\Bigl[\beta_k^2+\frac{2\beta_kM\delta_k}{|\D(x^{k-1},d^{k-1})|}\Bigr]\norm{d^{k-1}}^2 .
\]
Dividing by $\omega_k^2$ and using (i) and (iii), we obtain
\[
s_k\le\frac{\norm{v(x^k)}^2}{\omega_k^2}+\Bigl(1+\frac{2M\delta_k}{b\,\omega_k}\Bigr)s_{k-1}.
\]
By (ii), $\omega_k\ge a\norm{v(x^k)}^2\ge a\gamma^2$, so
$\norm{v(x^k)}^2/\omega_k^2\le1/(a\omega_k)\le1/(a^2\gamma^2)$ and
$2M\delta_k/(b\omega_k)\le\xi_k:=2M\delta_k/(ab\gamma^2)$. Therefore
\[s_k\le(1+\xi_k)s_{k-1}+\frac{1}{a^2\gamma^2} \quad \text{for}\quad k\ge1,\]
 and $s_0=\norm{v(x^0)}^2/\omega_0^2\le1/(a^2\gamma^2)$.
Moreover, $\sum_{j\ge1}\xi_j<\infty$ by Lemma~\ref{lem:step}(c). Let
$P_0:=1$ and $P_k:=\prod_{j=1}^k(1+\xi_j)$.
Since $\xi_j\ge0$ and
$1+\xi_j\le e^{\xi_j}$, we have $1\le P_k\le P:=\exp\bigl(\sum_{j\ge1}\xi_j\bigr)<\infty$.
Dividing the recurrence by $P_k$ gives
\[
\frac{s_k}{P_k}
\le \frac{s_{k-1}}{P_{k-1}}+\frac{1}{a^2\gamma^2P_k}
\le \frac{s_{k-1}}{P_{k-1}}+\frac{1}{a^2\gamma^2},
\]
since $P_k\ge1$. Summing these inequalities yields
\[
\frac{s_k}{P_k}\le s_0+\frac{k}{a^2\gamma^2}\le\frac{1}{a^2\gamma^2}(k+1).
\]
Consequently,
\[
s_k\le\frac{P}{a^2\gamma^2}(k+1).
\]
Then, by (iii),
\[
\sum_{k\ge0}\frac{\D(x^k,d^k)^2}{\norm{d^k}^2}\ \ge\ b^2\sum_{k\ge0}\frac1{s_k}
\ \ge\ \frac{a^2b^2\gamma^2}{P}\sum_{k\ge0}\frac1{k+1}=\infty ,
\]
which contradicts the Zoutendijk condition in Proposition~\ref{lem:zout}.
\end{proof}

We conclude this section by noting that conditions (ii)
and (iii) imply the sufficient descent condition
\eqref{eq:sd} with $\kappa=ab$.


\section{Convergence of the Fletcher--Reeves Method}\label{sec:fr}

As shown by Example~\ref{ex:fr}, the strong Wolfe conditions
do not guarantee descent for the FR method.
Sufficient descent can nevertheless be enforced by controlling
the tolerance in \eqref{eq:w2} relative to
$\norm{v(x^k)}^2$:
\begin{equation}\label{eq:wmu}
\varepsilon_k\le\mu\norm{v(x^k)}^2,
\qquad 0\le\mu<1.
\end{equation}
We also consider tolerances bounded by a fixed multiple of
$|\D(x^k,d^k)|$,
\begin{equation}\label{eq:gw}
\varepsilon_k\le-c_3\D(x^k,d^k),
\qquad c_3\ge0.
\end{equation}
With equality in \eqref{eq:gw}, \eqref{eq:w2} is the vector
counterpart of the generalized Wolfe conditions of
\cite{DY96} (see also \cite{HZ06}).
In particular, the choice
$\varepsilon_k=-c_2\D(x^k,d^k)$
recovers the strong Wolfe conditions.

Theorem~\ref{thm:fr} shows that \eqref{eq:wmu} yields sufficient descent and
global convergence. Under \eqref{eq:gw}, as in the scalar case \cite{DY96},
convergence holds assuming only that the search directions are $K$-descent
directions.

\begin{theorem}\label{thm:fr}
Assume that \textup{(A1)} and \textup{(A2)} hold, and consider
the method \eqref{eq:it}--\eqref{eq:dir} with
$0\le\beta_k\le\beta_k^{\FR}$ for all $k\ge1$.
Then
\[
\liminf_{k\to\infty}\norm{v(x^k)}=0
\]
in each of the following cases:
\begin{enumerate}[label=\textup{(\alph*)}]
\item The step sizes satisfy \eqref{eq:w1}--\eqref{eq:w2}
with tolerances satisfying \eqref{eq:wmu}. In this case, the method is well defined and
\[
\D(x^k,d^k)\le-(1-\mu)\norm{v(x^k)}^2
\qquad\text{for all }k\ge0.
\]
\item The search directions satisfy the sufficient descent
condition \eqref{eq:sd}, and the step sizes satisfy
\eqref{eq:w1}--\eqref{eq:w2} with arbitrary tolerances
$\varepsilon_k\in[0,+\infty]$.
\item Every search direction is a $K$-descent direction,
and the step sizes satisfy \eqref{eq:w1}--\eqref{eq:w2}
with tolerances satisfying \eqref{eq:gw}.
\end{enumerate}
These statements apply, in particular, to the pure FR
method, $\beta_k=\beta_k^{\FR}$.
\end{theorem}
\begin{proof}
We first prove (b), and then deduce (a) from it. Part (c) requires a separate
argument.

(b) With $\omega_k=\norm{v(x^k)}^2$, condition (i) of Theorem~\ref{thm:general}
reads $\beta_k\le\beta_k^{\FR}$, (ii) holds with $a=1$, and (iii) is the sufficient
descent condition \eqref{eq:sd} with $b=\kappa$. Hence the result follows from
Theorem~\ref{thm:general}.

(a) Since $x^0\in\mathcal L$ and
$\D(x^0,d^0)=-\norm{v(x^0)}^2$, the initial direction
is a $K$-descent direction satisfying the claimed estimate.
Assume inductively that $x^{k-1}\in\mathcal L$ and
$d^{k-1}$ is a $K$-descent direction.
Proposition~\ref{lem:exist} guarantees the existence of
admissible step sizes, and \eqref{eq:w1} implies
$x^k\in\mathcal L$.
By \eqref{eq:basic}, $0\le\beta_k\le\beta_k^{\FR}$, \eqref{eq:w2}, and \eqref{eq:wmu},
\[
\D(x^k,d^k)
\le-\norm{v(x^k)}^2
+\beta_k^{\FR}\max\{\D(x^k,d^{k-1}),0\}
\le-(1-\mu)\norm{v(x^k)}^2.
\]
Thus, by induction, the method is well defined and
\eqref{eq:sd} holds with $\kappa=1-\mu$.
The convergence conclusion follows from (b).

(c) Suppose, by contradiction, that
$\norm{v(x^k)}\ge\gamma>0$ for all $k$, and define
\[
r_k:=\frac{-\D(x^k,d^k)}{\norm{v(x^k)}^2},
\qquad
t_k:=\frac{\norm{d^k}^2}{\norm{v(x^k)}^4},
\qquad
\Lambda:=\frac{2M}{\gamma^2},
\]
where $M$ is the constant from Lemma~\ref{lem:key}.
By \eqref{eq:w1}, $x^k\in\mathcal L$ for all $k$, so Lemma~\ref{lem:key}
applies at every $k\ge1$.
Since every $d^k$ is a $K$-descent direction, $r_k>0$.
By \eqref{eq:basic},
$0\le\beta_k\le\beta_k^{\FR}$, \eqref{eq:w2}, and
\eqref{eq:gw} at iteration $k-1$,
\[
\begin{aligned}
\D(x^k,d^k)
&\le-\norm{v(x^k)}^2
+\beta_k^{\FR}\max\{\D(x^k,d^{k-1}),0\}\\
&\le-\norm{v(x^k)}^2(1-c_3r_{k-1}).
\end{aligned}
\]
Consequently,
\[
R_k:=r_k+c_3r_{k-1}\ge1
\qquad\text{for all }k\ge1.
\]

We next derive two estimates for $t_k$.
By \eqref{eq:dir},
\[
\norm{d^k}^2
=\norm{v(x^k)}^2
+2\beta_k\ip{v(x^k)}{d^{k-1}}
+\beta_k^2\norm{d^{k-1}}^2.
\]
Dividing by $\norm{v(x^k)}^4$ and using Lemma~\ref{lem:key} together with $0\le\beta_k\le\beta_k^{\FR}$ yields
\[
t_k
\le t_{k-1}
+\frac1{\norm{v(x^k)}^2}
+\frac{2M\delta_k}
       {r_{k-1}\norm{v(x^k)}^2}\,t_{k-1}.
\]
Multiplying by $r_{k-1}>0$ and using
$\norm{v(x^k)}\ge\gamma$ gives
\begin{equation}\label{eq:rk1tk}
r_{k-1}t_k
\le(r_{k-1}+\Lambda\delta_k)t_{k-1}
+\frac{r_{k-1}}{\gamma^2}.
\end{equation}

For the second estimate, \eqref{eq:dir} also gives
\[
\norm{d^k}^2
=2\ip{v(x^k)}{d^k}
-\norm{v(x^k)}^2
+\beta_k^2\norm{d^{k-1}}^2.
\]
Applying Lemma~\ref{lem:key} with $k+1$ in place of $k$ yields
\[
\ip{v(x^k)}{d^k}
\le M\delta_{k+1}
\frac{\norm{d^k}^2}{|\D(x^k,d^k)|}.
\]
Substituting this bound, dropping the negative term
$-\norm{v(x^k)}^2$, using
$\beta_k\le\beta_k^{\FR}$, and dividing by
$\norm{v(x^k)}^4$, we obtain
\[
t_k
\le t_{k-1}
+\frac{2M\delta_{k+1}}
       {r_k\norm{v(x^k)}^2}\,t_k.
\]
Multiplying by $r_k>0$ and using
$\norm{v(x^k)}\ge\gamma$ gives
\begin{equation}\label{eq:rktk}
r_kt_k
\le r_kt_{k-1}+\Lambda\delta_{k+1}t_k.
\end{equation}

Multiplying \eqref{eq:rk1tk} by $c_3$, adding \eqref{eq:rktk}, and using the
definition of $R_k$, we obtain
\[
R_kt_k
\le R_kt_{k-1}
+\frac{c_3r_{k-1}}{\gamma^2}
+c_3\Lambda\delta_k t_{k-1}
+\Lambda\delta_{k+1}t_k.
\]
Since $R_k\ge1$ and $c_3r_{k-1}\le R_k$, division
by $R_k$ gives
\[
(1-\Lambda\delta_{k+1})t_k
\le(1+c_3\Lambda\delta_k)t_{k-1}
+\frac1{\gamma^2}.
\]
By Lemma~\ref{lem:step}(c),
$\sum_k\delta_k<\infty$. Thus, there is $k_0$ such
that $\Lambda\delta_{k+1}\le1/2$ for all $k>k_0$.
For these indices,
\[
\frac{1+c_3\Lambda\delta_k}{1-\Lambda\delta_{k+1}}
=
1+\frac{\Lambda(c_3\delta_k+\delta_{k+1})}
        {1-\Lambda\delta_{k+1}}
\le1+2\Lambda(c_3\delta_k+\delta_{k+1}).
\]
Defining $\xi_k:=2\Lambda(c_3\delta_k+\delta_{k+1})$, it follows that
\[
t_k\le(1+\xi_k)t_{k-1}+\frac2{\gamma^2}.
\]
Since $\sum_k\xi_k<\infty$, the same product argument
as in the proof of Theorem~\ref{thm:general}
gives a constant $A>0$ such that
$t_k\le A(k+1)$ for all $k\ge k_0$.

Finally, $R_k\ge1$ implies
\[
\max\{r_{k-1},r_k\}\ge\frac1{1+c_3}.
\]
Thus, each of the disjoint pairs $\{2j,2j+1\}$
contains an index $k$ such that $r_k\ge1/(1+c_3)$.
Since
$\D(x^k,d^k)^2/\norm{d^k}^2=r_k^2/t_k$,
for any integer $J$ with $2J\ge k_0$ we obtain
\[
\sum_{k\ge2J}
\frac{\D(x^k,d^k)^2}{\norm{d^k}^2}
\ge
\frac1{A(1+c_3)^2}
\sum_{j\ge J}\frac1{2j+2}
=\infty,
\]
contradicting the Zoutendijk condition in Proposition~\ref{lem:zout}.
\end{proof}

Part (a) provides a concrete line search condition
under which the pure FR method is well defined
and globally convergent. Sufficient descent is
guaranteed by the line search and need not be
enforced separately.
Parts (b) and (c) cover the pure FR method in the two descent regimes considered
in \cite[Theorem~5.1]{LPP18}, sufficient descent and descent, where only strict
fractions of the FR parameter were allowed.

Part (c) is also the vector counterpart of
\cite[Theorem~2.3]{DY96}, see also \cite{DHLSYY00}.
Its proof follows the same consecutive-pair strategy.
The main difference lies in estimating the growth of
$\norm{d^k}^2/\norm{v(x^k)}^4$.
In the scalar case, this estimate follows from an exact identity for
$\D(x^k,d^k)/\norm{v(x^k)}^2$ and the curvature condition.
Here, it relies on the two estimates of
Lemma~\ref{lem:key} and the finite total decrease of $F$.

Example~\ref{ex:frcompact} in Section~\ref{sec:nonconvergence}
shows that, under the standard Wolfe conditions, descent alone
does not suffice, even with a compact initial level set.

\begin{remark}
(a) For the pure FR method in the scalar case,
a strong Wolfe line search with $c_2<1/2$
gives Al-Baali's bound
\[
|\D(x^k,d^k)|
\le\frac{\norm{v(x^k)}^2}{1-c_2}.
\]
Hence the strong Wolfe tolerance
$\varepsilon_k=c_2|\D(x^k,d^k)|$
satisfies \eqref{eq:wmu} with
$\mu=c_2/(1-c_2)<1$.
Thus Theorem~\ref{thm:fr}(a) recovers the classical
convergence result of \cite{AlB85}.

(b) Theorem~\ref{thm:fr} also applies to hybrid parameters
such as
\[
\beta_k=
\max\bigl\{0,\min\{\beta_k^{\PRP},\beta_k^{\FR}\}\bigr\},
\]
where $\beta_k^{\PRP}$ is the
Polak--Ribi\`ere--Polyak parameter of \cite{LPP18}.
Moreover, Theorem~\ref{thm:fr}(b) applies to parameters bounded by
$\beta_k^{\FR}$ whose sufficient descent comes from their construction, such as
the hybrids of \cite{YK24} and \cite{LCPO26} under suitable choices of their
parameters. In this case arbitrary tolerances are allowed in \eqref{eq:w2},
whereas both papers use the strong Wolfe conditions.

\end{remark}

\section{Convergence of the Dai--Yuan Method}\label{sec:dy}

In the scalar case, the DY method generates descent
directions and is globally convergent under the standard
Wolfe conditions \cite{DY99}.
In the vector setting, descent is still guaranteed, but, as
Example~\ref{ex:dy} in Section~\ref{sec:nonconvergence}
shows, global convergence may fail under the standard
Wolfe conditions. Theorem~\ref{thm:dy} shows that a tolerance in
\eqref{eq:w2} satisfying \eqref{eq:gw} ensures sufficient descent and global
convergence.

\begin{theorem}\label{thm:dy}
Assume that \textup{(A1)} and \textup{(A2)} hold, and consider
the method \eqref{eq:it}--\eqref{eq:dir} with
$0\le\beta_k\le\beta_k^{\DY}$ for all $k\ge1$.
Suppose that the step sizes satisfy
\eqref{eq:w1}--\eqref{eq:w2} with arbitrary tolerances
$\varepsilon_k\in[0,+\infty]$.
Then the method is well defined and every search direction
is a $K$-descent direction.
Moreover,
\[
\liminf_{k\to\infty}\norm{v(x^k)}=0
\]
in each of the following cases:
\begin{enumerate}[label=\textup{(\alph*)}]
\item The tolerances satisfy \eqref{eq:gw}.
In this case,
\[
\D(x^k,d^k)
\le-\frac{\norm{v(x^k)}^2}{1+c_3}
\qquad\text{for all }k\ge0.
\]
\item The search directions satisfy the sufficient descent
condition \eqref{eq:sd}.
\end{enumerate}
In particular, the pure DY method,
$\beta_k=\beta_k^{\DY}$, is globally convergent
under the strong Wolfe conditions.
\end{theorem}

\begin{proof}
We first show that the method is well defined and
generates $K$-descent directions.
Since $x^0\in\mathcal L$ and
$\D(x^0,d^0)=-\norm{v(x^0)}^2$, the initial
direction is a $K$-descent direction.
Let $k\ge1$ and assume that $x^{k-1}\in\mathcal L$
and $d^{k-1}$ is a $K$-descent direction.
Proposition~\ref{lem:exist} guarantees the existence of
admissible step sizes, and \eqref{eq:w1} implies
$x^k\in\mathcal L$.
By \eqref{eq:w2},
\[
\D(x^k,d^{k-1})-\D(x^{k-1},d^{k-1})
\ge(1-c_2)|\D(x^{k-1},d^{k-1})|>0,
\]
so $\beta_k^{\DY}>0$ is well defined.
Since $0\le\beta_k\le\beta_k^{\DY}$,
\[
\norm{v(x^k)}^2
\ge
\beta_k\bigl(
\D(x^k,d^{k-1})-\D(x^{k-1},d^{k-1})
\bigr).
\]
Combining this inequality with \eqref{eq:basic} gives
\[
\D(x^k,d^k)
\le-\norm{v(x^k)}^2
+\beta_k\D(x^k,d^{k-1})
\le\beta_k\D(x^{k-1},d^{k-1}).
\]
If $\beta_k>0$, the last expression is negative.
If $\beta_k=0$, then $d^k=v(x^k)$ and
$\D(x^k,d^k)=-\norm{v(x^k)}^2<0$.
Thus, by induction, the method is well defined and
every search direction is a $K$-descent direction.
The preceding inequality also yields
\[
\beta_k
\le
\frac{\D(x^k,d^k)}
     {\D(x^{k-1},d^{k-1})}
\qquad\text{for all }k\ge1,
\]
because $\D(x^{k-1},d^{k-1})<0$.
We now prove (b), and then deduce (a).

(b) With $\omega_k=-\D(x^k,d^k)$, condition (i) of
Theorem~\ref{thm:general} is the ratio bound above,
(ii) holds with $a=\kappa$ by \eqref{eq:sd}, and
(iii) holds with $b=1$. Hence the result follows from
Theorem~\ref{thm:general}.

(a) The claimed sufficient descent estimate holds
for $k=0$.
Let $k\ge1$.
If $\D(x^k,d^{k-1})\le0$, \eqref{eq:basic} gives
\[
\D(x^k,d^k)\le-\norm{v(x^k)}^2
\le-\frac{\norm{v(x^k)}^2}{1+c_3}.
\]
If $\D(x^k,d^{k-1})>0$, then
\eqref{eq:w2} and \eqref{eq:gw} at iteration $k-1$
give
$\D(x^k,d^{k-1})
\le c_3|\D(x^{k-1},d^{k-1})|$.
Hence, by \eqref{eq:basic} and
$\beta_k\le\beta_k^{\DY}$,
\[
\begin{aligned}
\D(x^k,d^k)
&\le-\norm{v(x^k)}^2
+\beta_k^{\DY}\D(x^k,d^{k-1})\\
&=-\norm{v(x^k)}^2
\frac{|\D(x^{k-1},d^{k-1})|}
     {\D(x^k,d^{k-1})
      +|\D(x^{k-1},d^{k-1})|}\\
&\le-\frac{\norm{v(x^k)}^2}{1+c_3}.
\end{aligned}
\]
Thus \eqref{eq:sd} holds with $\kappa=1/(1+c_3)$,
and convergence follows from (b).
\end{proof}

For $c_3=c_2$, the sufficient descent estimate in (a) is
\cite[Lemma~5.5]{LPP18}, and the ratio bound in the proof follows as in the
proof of \cite[Theorem~5.7]{LPP18}. Convergence under the standard Wolfe
conditions and sufficient descent for parameters bounded by
$\delta\,\D(x^k,d^k)/\D(x^{k-1},d^{k-1})$ with $0\le\delta<1$ was established in
\cite[Theorem~3]{HLZC24}.

This criterion also implies convergence for
$\beta_k=\eta\beta_k^{\DY}$ with $0\le\eta<1$, under the assumptions of that
paper. Indeed, for $0<\eta<1$, \eqref{eq:basic} and the definition of
$\beta_k^{\DY}$ give
\[
\D(x^k,d^k)
\le-(1-\eta)\norm{v(x^k)}^2
+\beta_k\D(x^{k-1},d^{k-1}),
\]
so \eqref{eq:sd} holds with $\kappa=1-\eta$.
Moreover, by the lower bound in \eqref{eq:w2},
\[
\begin{aligned}
\frac{\norm{v(x^k)}^2}{|\D(x^{k-1},d^{k-1})|}
&=\frac{\beta_k}{\eta}
\frac{\D(x^k,d^{k-1})-\D(x^{k-1},d^{k-1})}
     {|\D(x^{k-1},d^{k-1})|}\\
&\ge\frac{1-c_2}{\eta}\,\beta_k.
\end{aligned}
\]
Since $\D(x^{k-1},d^{k-1})<0$, these two inequalities imply
\[
\beta_k\le\delta_\eta\,
\frac{\D(x^k,d^k)}{\D(x^{k-1},d^{k-1})},
\qquad
\delta_\eta:=\frac{\eta}{1-c_2(1-\eta)}<1.
\]
Thus \cite[Theorem~3]{HLZC24} applies. The case $\eta=0$ is the steepest descent
method. For the pure DY method, $\eta=1$ gives $\delta_\eta=1$.
Theorem~\ref{thm:general}, applied with $\omega_k=-\D(x^k,d^k)$, allows
$\delta=1$, and the new ingredient is Lemma~\ref{lem:key}.

Theorem~\ref{thm:dy} covers several parameters at once. Recall that the
conjugate descent (CD) parameter of \cite{LPP18} is
\[
\beta_k^{\CD}=\frac{\norm{v(x^k)}^2}{-\D(x^{k-1},d^{k-1})}.
\]

\begin{corollary}\label{cor:dy}
Assume that \textup{(A1)} and \textup{(A2)} hold, and that the step sizes satisfy
\eqref{eq:w1}--\eqref{eq:w2} with tolerances satisfying \eqref{eq:gw}. Then
$\liminf_{k\to\infty}\norm{v(x^k)}=0$ for each of the following choices:
\begin{enumerate}[label=\textup{(\alph*)}]
\item $\beta_k=\eta\beta_k^{\DY}$ with $0\le\eta\le1$;
\item the hybrid parameter
\[
\beta_k=\max\bigl\{0,\min\{\beta_k^{\HS},\beta_k^{\DY}\}\bigr\},
\]
where
 $\beta_k^{\HS}$ is the Hestenes--Stiefel parameter of \cite{LPP18};
\item the hybrid parameter of \cite{PZY26}, equal to $\beta_k^{\DY}$ if
 $\D(x^k,d^{k-1})>0$ and to $\beta_k^{\CD}$ otherwise, that is,
 $\beta_k=\min\{\beta_k^{\CD},\beta_k^{\DY}\}$;
\item the Dai--Yuan vector family $\beta_k(\tau_k)$ of \cite{HLZC24}, for every
 choice of $\tau_k\in\R$ and of the constant $\mu\ge1$ in \cite[(25)]{HLZC24};
\item $\beta_k=\eta\beta_k^{\CD}$ with $0\le\eta\le1/(1+c_3)$.
\end{enumerate}
\end{corollary}
\begin{proof}
In all cases $0\le\beta_k\le\beta_k^{\DY}$, so the result follows from
Theorem~\ref{thm:dy}(a). This is clear for (a) and (b). For (c), note that if
$\D(x^k,d^{k-1})\le0$, then the denominator of $\beta_k^{\DY}$ is at most
$-\D(x^{k-1},d^{k-1})$, so $\beta_k^{\CD}\le\beta_k^{\DY}$. For (d), \cite[Theorem~2 and (21)]{HLZC24} gives
\[
0<\beta_k(\tau_k)
\le
\frac{\norm{v(x^k)}^2}
{\D(x^k,d^{k-1})
+\max\{\D(x^k,-d^{k-1}),0\}
-\D(x^{k-1},d^{k-1})}
\le\beta_k^{\DY}.
\] 
For (e), \eqref{eq:w2} and \eqref{eq:gw} give
$\D(x^k,d^{k-1})\le c_3|\D(x^{k-1},d^{k-1})|$, so the denominator of
$\beta_k^{\DY}$ is at most $(1+c_3)|\D(x^{k-1},d^{k-1})|$ and
$\beta_k^{\CD}\le(1+c_3)\beta_k^{\DY}$.
\end{proof}

Under the strong Wolfe conditions, item (e)
allows $0\le\eta\le1/(1+c_2)$, which contains the range $0\le\eta<1-c_2$ of
\cite[Theorem~5.4]{LPP18}. 
Taking $c_3=0$, so that $\varepsilon_k=0$,
also covers the pure CD method.
This parallels the scalar case, where the CD method converges under the generalized
Wolfe conditions with zero upper tolerance and may fail to converge under the strong
Wolfe conditions (see \cite{HZ06} and the references therein).

Theorem~\ref{thm:dy}(b) also applies to parameters bounded by $\beta_k^{\DY}$ whose
sufficient descent comes from their construction. For instance, the parameter of
\cite{PCLL25} satisfies $0<\beta_k\le\mu_1\beta_k^{\DY}$ and, with
$\mu_2=\mu_1/c_2$, yields \eqref{eq:sd} with $\kappa=1-c_2$ under the lower bound in \eqref{eq:w2}
\cite[Lemma~3.2]{PCLL25}. Hence Theorem~\ref{thm:dy}(b) gives convergence under the
standard Wolfe conditions for every $\mu_1\in(0,1]$, whereas
\cite[Theorem~3.1]{PCLL25} requires $\mu_1<\sqrt{2/3}\,(1-c_2)^2$.

\section{Failure of Convergence under the Standard Wolfe Conditions}\label{sec:nonconvergence}

We now show that the standard Wolfe conditions do not ensure
global convergence of the pure FR and DY methods, even when
every search direction is a $K$-descent direction and the initial
level set is compact. Both examples are smooth biobjective
problems in dimension two and use the same sequence of points.

Throughout this section, let $K=\R^2_+$, take $C$ as the
canonical basis of $\R^2$, and set $e=(1,1)^\top$.

\subsection{A Common Construction}\label{sec:counterconstruction}

We first prescribe the sequence of points and then construct
the objective functions around it. The local expressions of
these functions will be chosen separately for FR and DY.

Let
\[
R:=\begin{pmatrix}0&-1\\1&0\end{pmatrix},
\qquad e_1:=(1,0)^\top.
\]
The matrix $R$ represents the counterclockwise rotation by $90^\circ$.
For any $u,v\in\R^2$ and $k\ge0$, we have
\begin{equation}\label{eq:counterrotation}
\ip{R^{k+1}u}{R^kv}=\ip{Ru}{v},
\qquad
\ip{Ru}{e_1}=-u_2.
\end{equation}

\begin{figure}[tbp]
\definecolor{cgblue}{HTML}{28577C}
\definecolor{cgteal}{HTML}{16857A}
\definecolor{cgorange}{HTML}{C77723}
\definecolor{cggray}{HTML}{87929B}
\tikzset{
 cg axis/.style={cggray!55,line width=.35pt,-{Stealth[length=1.5mm]}},
 cg path/.style={cgblue,line width=.85pt},
 cg flow/.style={cg path,postaction={decorate},
   decoration={markings,mark=at position .58 with
   {\arrow{Stealth[length=1.6mm,width=1.2mm]}}}},
 cg label/.style={font=\small,inner sep=2pt},
 cg tiny/.style={font=\footnotesize,inner sep=2pt}
}
\centering
\resizebox{0.85\linewidth}{!}{%
\begin{minipage}{\textwidth}
\begin{minipage}[t]{.55\textwidth}
\centering
\textbf{\small (a) The common trajectory}\par\medskip
\begin{tikzpicture}[x=2.75cm,y=2.75cm]
 \path[use as bounding box] (-1.12,-.97) rectangle (1.02,.94);
 \foreach \j/\xx/\yy in {0/-.5/-.5,1/.5/-.5,2/.5/.5,3/-.5/.5}{
   \coordinate (Y\j) at (\xx,\yy);
   \fill[cgblue!5] (Y\j) circle[radius=.25];
   \draw[cgblue!20,line width=.35pt] (Y\j) circle[radius=.25];
 }
 \draw[cg axis] (-1.02,0)--(.94,0) node[cg tiny,below] {$x_1$};
 \draw[cg axis] (0,-.84)--(0,.84) node[cg tiny,left] {$x_2$};
 \node[cg tiny,below left,cggray] at (0,0) {$0$};
 \draw[cggray,densely dashed,line width=.65pt]
   (Y0)--(Y1)--(Y2)--(Y3)--cycle;
 \foreach \k in {0,...,13}{
   \pgfmathsetmacro{\ee}{.1*pow(2,-\k)}
   \pgfmathtruncatemacro{\j}{mod(\k,4)}
   \ifcase\j
     \pgfmathsetmacro{\xx}{-.5-4*\ee}
     \pgfmathsetmacro{\yy}{-.5-2*\ee}
   \or
     \pgfmathsetmacro{\xx}{.5+2*\ee}
     \pgfmathsetmacro{\yy}{-.5-4*\ee}
   \or
     \pgfmathsetmacro{\xx}{.5+4*\ee}
     \pgfmathsetmacro{\yy}{.5+2*\ee}
   \or
     \pgfmathsetmacro{\xx}{-.5-2*\ee}
     \pgfmathsetmacro{\yy}{.5+4*\ee}
   \fi
   \coordinate (X\k) at (\xx,\yy);
 }
 \foreach \k in {0,...,6}{
   \pgfmathtruncatemacro{\next}{\k+1}
   \draw[cg flow] (X\k)--(X\next);
 }
 \foreach \k in {7,...,12}{
   \pgfmathtruncatemacro{\next}{\k+1}
   \draw[cg path,opacity=.45,line width=.55pt] (X\k)--(X\next);
 }
 \foreach \k in {0,...,7}{\fill[cgblue] (X\k) circle[radius=.019];}
 \foreach \j in {0,...,3}{
   \fill[white] (Y\j) circle[radius=.022];
   \draw[cgteal,line width=.9pt] (Y\j) circle[radius=.022];
 }
 \node[cg label,above right,cgteal] at (Y0) {$y_0$};
 \node[cg label,above left,cgteal] at (Y1) {$y_1$};
 \node[cg label,below left,cgteal] at (Y2) {$y_2$};
 \node[cg label,below right,cgteal] at (Y3) {$y_3$};
 \node[cg label,below left] at (X0) {$x^0$};
 \node[cg label,below right] at (X1) {$x^1$};
 \node[cg label,above right] at (X2) {$x^2$};
 \node[cg label,above left] at (X3) {$x^3$};
 \node[cg tiny,cggray] at (0,-.94) {Accumulation points: the four vertices};
\end{tikzpicture}
\end{minipage}\hfill
\begin{minipage}[t]{.42\textwidth}
\centering
\textbf{\small (b) Iterates near $y_0$}\par\medskip
\begin{tikzpicture}[x=3.15cm,y=3.15cm]
 \path[use as bounding box] (-1.18,-.85) rectangle (.56,.82);
 \coordinate (L) at (0,0);
 \coordinate (A) at (-1,-.5);
 \coordinate (B) at (-.0625,-.03125);
 \coordinate (C) at (-.00390625,-.001953125);
 \draw[cggray,densely dashed,line width=.65pt] (-1.14,0)--(.51,0);
 \draw[cggray,densely dashed,line width=.65pt] (0,-.71)--(0,.53);
 \draw[cgorange,densely dashed,line width=.8pt] (-1.09,-.545)--(L);
 \draw[cg flow] (-1,.48)--(A);
 \draw[cg flow] (A)--(.51,-.5);
 \draw[cg path,opacity=.6,line width=.65pt] (-.0625,.48)--(B)--(.51,-.03125);
 \draw[cg path,opacity=.35,line width=.5pt] (-.00390625,.48)--(C)--(.51,-.001953125);
 \fill[cgblue] (A) circle[radius=.016];
 \fill[cgblue] (B) circle[radius=.013];
 \fill[white] (L) circle[radius=.021];
 \draw[cgteal,line width=.95pt] (L) circle[radius=.021];
 \node[cg label,below left] at (A) {$x^4$};
 \draw[cggray,line width=.45pt] (B)--(-.38,.27);
 \node[cg label,above left,fill=white,inner sep=1pt] at (-.38,.27) {$x^8$};
 \node[cg label,above right,cgteal,fill=white,inner sep=1pt] at (L) {$y_0$};
 \node[cg tiny,cgorange,below right,rotate=26.565] at (-.72,-.36)
   {approach along $p$};
 \node[cg tiny,align=center] at (-.31,-.79)
   {$x^{k+4}-y_0=\tfrac1{16}(x^k-y_0)$\\[2pt]$k\equiv0\pmod4$};
\end{tikzpicture}
\end{minipage}
\end{minipage}%
}
\caption{The sequence of points shared by the FR and DY examples.
(a) The first steps and the four accumulation points. Each shaded disk has radius
$1/4$, where the corresponding cutoff equals one. The dashed square joins the
four accumulation points. (b) Successive turns near $y_0$: after every four steps,
the distance to the same vertex is divided by $16$.}
\label{fig:countertrajectory}
\end{figure}
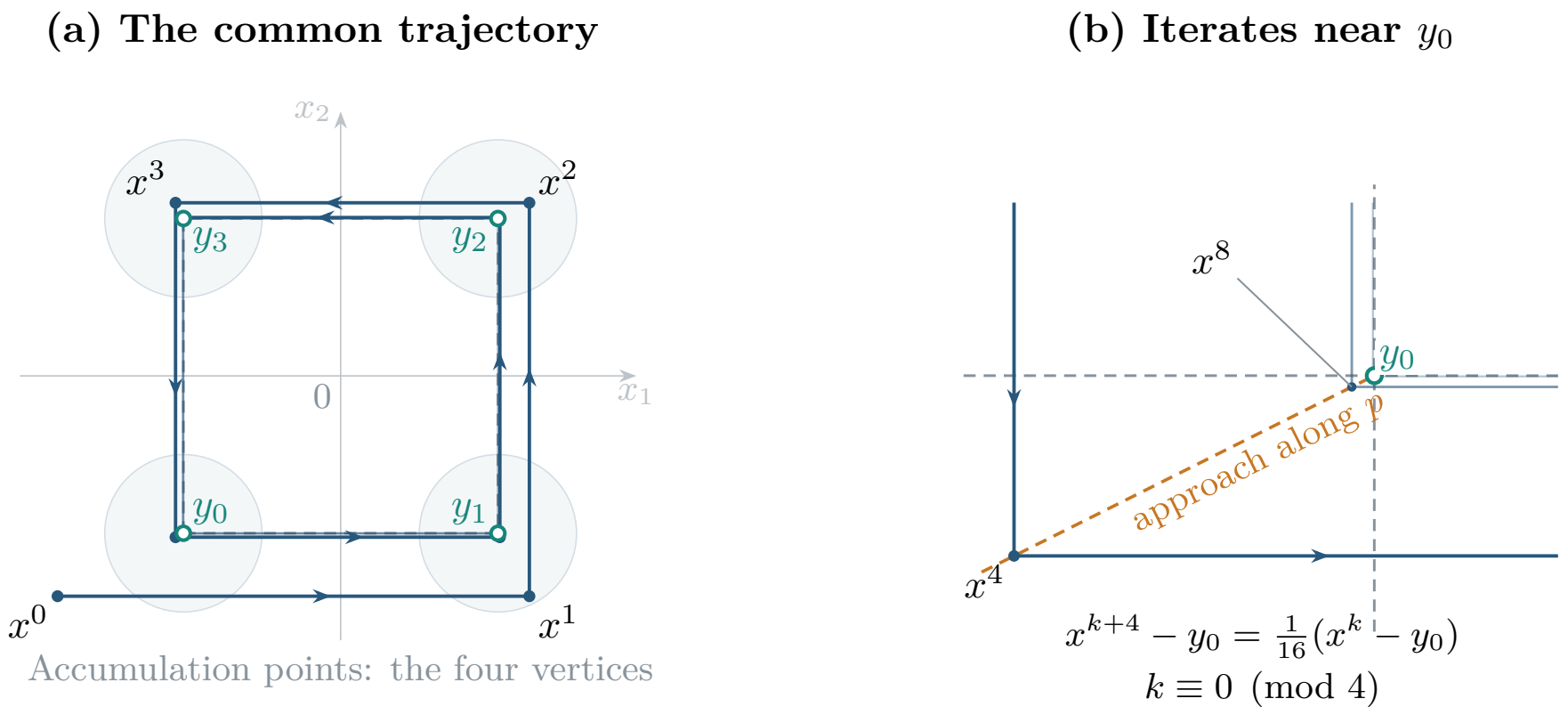

To define the trajectory, set
\[
y_0:=(-\tfrac12,-\tfrac12)^\top,
\qquad y_j:=R^jy_0\quad(j=1,2,3),
\qquad p:=(-4,-2)^\top.
\]
Then define
\begin{equation}\label{eq:countertrajectory}
\zeta_k:=\frac{2^{-k}}{10},
\qquad
x^k:=R^k(y_0+\zeta_kp),
\qquad k\ge0.
\end{equation}
Since $Ry_0-y_0=e_1$ and $\tfrac12Rp-p=5e_1$,
\begin{equation}\label{eq:countergeometry}
x^{k+1}-x^k=(1+5\zeta_k)R^ke_1.
\end{equation}
Thus consecutive displacements follow the four coordinate
directions in cyclic order. In particular,
\[
x^0=(-9/10,-7/10)^\top,
\qquad
x^1=(3/5,-7/10)^\top.
\]
For each $j\in\{0,1,2,3\}$, the subsequence
$\{x^{4\ell+j}\}$ converges to $y_j$.
These four vertices are the only accumulation points
of the sequence, as illustrated in Figure~\ref{fig:countertrajectory}.

To define the objective functions, let $B(y,r)$ denote
the open Euclidean ball with center $y$ and radius $r$.
Choose $C^\infty$ functions
$\chi_0,\chi_1,\chi_2,\chi_3\colon\R^2\to[0,1]$ such that,
for $j=0,1,2,3$,
\[
\chi_j(x)=1 \quad\text{if }\norm{x-y_j}\le1/4,
\qquad
\chi_j(x)=0 \quad\text{if }\norm{x-y_j}\ge3/10.
\]
For the initial point $x^0$, choose a $C^\infty$ function
$\chi_{\mathrm{init}}\colon\R^2\to[0,1]$ that equals one
near $x^0$ and vanishes whenever $\norm{x-x^0}\ge1/20$.
The balls $B(y_j,3/10)$ and $B(x^0,1/20)$
are pairwise disjoint. Indeed, distinct vertices are
at least one unit apart, and
\[
\min_{0\le j\le3}\norm{x^0-y_j}
=\norm{x^0-y_0}
=\frac{\sqrt{20}}{10}
>\frac3{10}+\frac1{20}.
\]
Consequently, at most one of the cutoff functions
can be nonzero at any point.

Let $h_1,h_2\colon\R^2\to\R$ be polynomials and let $\gamma,\tau>0$.
Define $F=(f_1,f_2)^\top\colon\R^2\to\R^2$ by
\begin{equation}\label{eq:counterextension}
\begin{aligned}
f_i(x):={}&\norm{x}^2
+\sum_{j=0}^3\chi_j(x)\bigl[h_i(R^{-j}x-y_0)-\norm{x}^2\bigr]\\
&+\chi_{\mathrm{init}}(x)
  \bigl[\gamma-\tau\ip{e_1}{x-x^0}-\norm{x}^2\bigr],
\qquad i=1,2.
\end{aligned}
\end{equation}
Thus each component agrees with the corresponding local
expression near the prescribed points and with $\norm{x}^2$
outside a bounded region.

\begin{lemma}\label{lem:counterextension}
Let $F$ be defined by \eqref{eq:counterextension}.
Then $F$ is $C^\infty$, its Jacobian is globally
Lipschitz continuous, and both components are coercive.
In particular, the initial level set $\mathcal L$ is compact,
and \textup{(A1)} and \textup{(A2)} hold.
Moreover, for $i=1,2$,
\begin{align}
f_i(x^0)&=\gamma,
&
\nabla f_i(x^0)&=-\tau e_1,
\label{eq:counterstart}\\
f_i(x^k)&=h_i(\zeta_kp),
&
\nabla f_i(x^k)&=R^k\nabla h_i(\zeta_kp),
\qquad k\ge1.
\label{eq:countervalues}
\end{align}
\end{lemma}

\begin{proof}
The expression in \eqref{eq:counterextension} is a finite
sum of smooth functions, so $F$ is $C^\infty$.
Outside a fixed bounded region, both components equal $\norm{x}^2$,
so each has Hessian $2I$, where $I$ is the $2\times2$ identity matrix.
On the closure of that region, their Hessians are bounded
by continuity and compactness. Hence the Hessians are
globally bounded, which implies that $JF$ is globally
Lipschitz continuous.
Both components are also coercive, so $\mathcal L$
is closed and bounded, and therefore compact.
Global Lipschitz continuity gives (A1), and continuity
of $F$ on the compact set $\mathcal L$ gives a componentwise
lower bound, so (A2) holds.

Near $x^0$, we have $\chi_{\mathrm{init}}=1$ and $\chi_j=0$
for $j=0,1,2,3$.
Thus $f_i(x)=\gamma-\tau\ip{e_1}{x-x^0}$ there,
which proves \eqref{eq:counterstart}.
For $k\ge1$,
\[
\norm{\zeta_kp}\le\frac{\sqrt{20}}{20}<\frac14.
\]
Let $j\in\{0,1,2,3\}$ satisfy $j\equiv k\pmod4$.
Then $x^k$ lies in $B(y_j,1/4)$.
Near this point, $\chi_j=1$ and all the other cutoff
functions vanish. Hence
\[
f_i(x)=h_i(R^{-k}x-y_0).
\]
Evaluating this expression and its gradient at $x^k$,
and using $(R^{-k})^\top=R^k$, proves
\eqref{eq:countervalues}.
\end{proof}

\subsection{The Fletcher--Reeves Method}\label{sec:frnonconvergence}

In Example~\ref{ex:fr}, the pure FR method can lose descent under
a strong Wolfe line search. We now show that even when every
direction is a $K$-descent direction and the initial level set
is compact, the standard Wolfe conditions do not ensure global
convergence.

\begin{example}\label{ex:frcompact}
In \eqref{eq:counterextension}, take $\gamma=4$, $\tau=2$, and
\[
h_1(s,t):=-2t+\tfrac12st-\tfrac34t^2,
\qquad
h_2(s,t):=-2s.
\]
By Lemma~\ref{lem:counterextension}, $JF$ is globally Lipschitz
continuous and $\mathcal L$ is compact. We show that the pure FR
method generates the sequence \eqref{eq:countertrajectory} with
the step sizes
\[
\alpha_0:=\tfrac34,
\qquad
\alpha_k:=1+5\zeta_k \quad(k\ge1).
\]

We first determine the steepest descent directions at the prescribed
points. At $x^0$, both gradients equal $-2e_1$, so $v(x^0)=2e_1$.
For $k\ge1$, the gradients at $x^k$ are obtained by rotating the
gradients of the local model at $\zeta_kp$, by \eqref{eq:countervalues}.
For $0<\zeta\le1/20$, these local gradients are
\[
g_1:=\nabla h_1(\zeta p)=(-\zeta,-2+\zeta)^\top,
\qquad
g_2:=\nabla h_2(\zeta p)=(-2,0)^\top.
\]
Since $g_1-g_2=(2-\zeta)(1,-1)^\top$, the gradient segment
is perpendicular to $e=(1,1)^\top$. Also,
\[
-e=g_2+\frac{1}{2-\zeta}(g_1-g_2).
\]
Since $0<1/(2-\zeta)<1$, this identity shows that $-e\in[g_1,g_2]$.
As $e$ is perpendicular to $g_1-g_2$, $-e$ is the projection
of the origin onto $[g_1,g_2]$.
By \eqref{eq:countervalues}, the gradient segment at $x^k$
is the rotation by $R^k$ of this local segment with $\zeta=\zeta_k$.
Its minimum-norm point is consequently $-R^ke$, so
\[
v(x^k)=R^ke,
\qquad
\norm{v(x^k)}^2=2
\qquad(k\ge1).
\]

We next show by induction that the prescribed step sizes satisfy
the standard Wolfe conditions and that the pure FR method generates
the points $x^k$, with $d^k=R^ke_1$ for all $k\ge1$.

The initial direction is $d^0=2e_1$, and
$x^0+\alpha_0d^0=x^1$ by \eqref{eq:countergeometry}.
For the first step,
\[
F(x^0)=(4,4)^\top,\qquad
F(x^1)=(81/400,2/5)^\top,\qquad
\alpha_0\D(x^0,d^0)=-3.
\]
Both objectives decrease by more than $3$, so \eqref{eq:w1}
holds for every $c_1\in(0,1)$. The curvature condition also
holds, since $\D(x^0,d^0)=-4$ and
$\D(x^1,d^0)=2(2-\zeta_1)=39/10>0$.
Thus $\alpha_0$ is admissible. Since $\beta_1^{\FR}=2/4=1/2$,
the FR update gives
\[
d^1=Re+\tfrac12(2e_1)=Re_1.
\]

Now suppose that the method has reached $x^k$ with $k\ge1$ and
$d^k=R^ke_1$. By \eqref{eq:countergeometry}, the candidate step
$\alpha_k$ reaches $x^{k+1}$. The local gradients and
\eqref{eq:counterrotation} give
\begin{equation}\label{eq:frcompactderivatives}
\D(x^k,d^k)=\max\{-\zeta_k,-2\}=-\zeta_k<0,
\qquad
\D(x^{k+1},d^k)=2-\zeta_{k+1}>0.
\end{equation}
Thus $d^k$ is a $K$-descent direction and the standard curvature
condition holds for every $c_2\in(0,1)$.

To verify the Armijo condition, observe that
\[
h_1(\zeta p)=4\zeta+\zeta^2,
\qquad
h_2(\zeta p)=8\zeta.
\]
Since $\zeta_{k+1}=\zeta_k/2$, the decreases in $f_1$ and $f_2$
are $2\zeta_k+3\zeta_k^2/4$ and $4\zeta_k$, respectively.
Both are at least $2\zeta_k$, whereas
\[
\alpha_k|\D(x^k,d^k)|
=(1+5\zeta_k)\zeta_k<2\zeta_k,
\]
because $\zeta_k\le1/20$. Hence \eqref{eq:w1} holds
for every $c_1\in(0,1)$.

Having verified that $\alpha_k$ is admissible, we compute the next
direction. Since $\norm{v(x^{k+1})}^2=\norm{v(x^k)}^2=2$,
we have $\beta_{k+1}^{\FR}=1$. Using $Re+e_1=Re_1$, we obtain
\[
d^{k+1}=R^{k+1}e+R^ke_1=R^{k+1}e_1.
\]
This completes the induction.

The same step sizes therefore satisfy the standard Wolfe
conditions for every $0<c_1<c_2<1$. Every search direction
is a $K$-descent direction, yet $\norm{v(x^k)}^2=2$ for
every $k\ge1$, so the method is not globally convergent.\hfill$\diamond$
\end{example}

The example also shows how the hypotheses of
Theorem~\ref{thm:fr}(b) and (c) fail. By \eqref{eq:frcompactderivatives},
\[
\frac{-\D(x^k,d^k)}{\norm{v(x^k)}^2}
=\frac{\zeta_k}{2}\longrightarrow0,
\qquad
\frac{\D(x^{k+1},d^k)}{|\D(x^k,d^k)|}
=\frac{2-\zeta_{k+1}}{\zeta_k}\longrightarrow\infty.
\]
The first limit rules out sufficient descent. The second shows
that no fixed $c_3$ can satisfy the upper curvature bound
along this sequence.

\subsection{The Dai--Yuan Method}\label{sec:dynonconvergence}

We now choose a local model for which the pure DY method follows
the same trajectory as the FR method in Example~\ref{ex:frcompact}.
The model makes the denominator of the DY parameter constant
after the first step.

\begin{example}\label{ex:dy}
In \eqref{eq:counterextension}, take $\gamma=1$, $\tau=1/10$, and
\[
h_1(s,t):=-\tfrac52t+\tfrac12st-t^2,
\qquad
h_2(s,t):=-\tfrac54s.
\]
By Lemma~\ref{lem:counterextension}, $JF$ is globally Lipschitz
continuous and $\mathcal L$ is compact. We use the step sizes
\[
\alpha_0:=15,
\qquad
\alpha_k:=1+5\zeta_k \quad(k\ge1).
\]

We first compute the steepest descent directions.
At $x^0$, both gradients equal $-e_1/10$, so $v(x^0)=e_1/10$.
For $0<\zeta\le1/20$, the local models have the same value
$h_1(\zeta p)=h_2(\zeta p)=5\zeta$, and their gradients are
\[
g_1:=\nabla h_1(\zeta p)=(-\zeta,-\tfrac52+2\zeta)^\top,
\qquad
g_2:=\nabla h_2(\zeta p)=(-\tfrac54,0)^\top.
\]
Let $u:=(1,1/2)^\top$. Since
$g_1-g_2=(5/4-\zeta)(1,-2)^\top$, the segment $[g_1,g_2]$
is perpendicular to $u$. Moreover,
\[
-u=g_2+\frac{1}{5-4\zeta}(g_1-g_2).
\]
The coefficient lies in $(0,1)$, so $-u$ belongs to $[g_1,g_2]$
and is the projection of the origin onto this segment.
Rotating the gradients as in \eqref{eq:countervalues} gives
\[
v(x^k)=R^ku,
\qquad
\norm{v(x^k)}^2=\tfrac54
\qquad(k\ge1).
\]

We next show by induction that the step sizes satisfy the standard
Wolfe conditions and that the pure DY method generates the points
in \eqref{eq:countertrajectory}, with $d^k=R^ke_1$ for $k\ge1$.
For the first step, $d^0=e_1/10$ and
$x^0+\alpha_0d^0=x^1$ by \eqref{eq:countergeometry}. Since
\[
F(x^0)=e,\qquad F(x^1)=\tfrac14e,\qquad
\alpha_0\D(x^0,d^0)=-\tfrac3{20},
\]
each objective decreases by $3/4$, so the Armijo condition
holds for every $c_1\in(0,1)$. The curvature condition also
holds, since
\[
\D(x^0,d^0)=-\tfrac1{100},
\qquad
\D(x^1,d^0)=\tfrac1{10}\bigl(\tfrac52-2\zeta_1\bigr)
=\tfrac6{25}>0.
\]
Thus $\alpha_0$ is admissible. The first DY update gives
\[
\beta_1^{\DY}=\frac{5/4}{6/25+1/100}=5,
\qquad
d^1=Ru+5(e_1/10)=Re_1.
\]

Now suppose that the method has reached $x^k$ with $k\ge1$ and
$d^k=R^ke_1$. By \eqref{eq:countergeometry}, the candidate step
$\alpha_k$ reaches $x^{k+1}$. Using \eqref{eq:counterrotation}
and $\zeta_{k+1}=\zeta_k/2$, we obtain
\begin{equation}\label{eq:dycompactderivatives}
\D(x^k,d^k)=\max\{-\zeta_k,-\tfrac54\}=-\zeta_k<0,
\qquad
\D(x^{k+1},d^k)=\tfrac52-\zeta_k>0.
\end{equation}
Thus $d^k$ is a $K$-descent direction and the standard curvature
condition holds for every $c_2\in(0,1)$.

For the Armijo condition, recall that $f_i(x^k)=5\zeta_k$
for $i=1,2$. Each objective therefore decreases by
$5(\zeta_k-\zeta_{k+1})=5\zeta_k/2$. Since $\zeta_k\le1/20$,
\[
\alpha_k|\D(x^k,d^k)|
=(1+5\zeta_k)\zeta_k
\le\tfrac54\zeta_k
<\tfrac52\zeta_k.
\]
Hence \eqref{eq:w1} holds for every $c_1\in(0,1)$.

By \eqref{eq:dycompactderivatives}, the denominator of the DY
parameter is $(5/2-\zeta_k)-(-\zeta_k)=5/2$. Consequently,
\[
\beta_{k+1}^{\DY}=\frac{5/4}{5/2}=\tfrac12,
\qquad
d^{k+1}=R^{k+1}u+\tfrac12R^ke_1=R^{k+1}e_1,
\]
where the last equality follows from $Ru+\tfrac12e_1=Re_1$.
This completes the induction.

The same step sizes satisfy the standard Wolfe conditions for
every $0<c_1<c_2<1$, and every search direction is a
$K$-descent direction. Nevertheless, $\norm{v(x^k)}^2=5/4$
for all $k\ge1$, so the pure DY method is not globally
convergent.\hfill$\diamond$
\end{example}

Examples~\ref{ex:frcompact} and~\ref{ex:dy} illustrate the
same gap in \eqref{eq:vd}. For every $k\ge1$, both examples give
\[
\ip{v(x^k)}{d^k}=1,\qquad
-\D(x^k,d^k)=\zeta_k\longrightarrow0.
\]
Since $\norm{d^k}=1$, the terms of the Zoutendijk series are
$\zeta_k^2$ for $k\ge1$, so the series converges even though
$\norm{v(x^k)}$ stays bounded away from zero.
For DY, we also have
\[
\frac{-\D(x^k,d^k)}{\norm{v(x^k)}^2}
=\frac{4\zeta_k}{5}\longrightarrow0,
\qquad
\frac{\D(x^{k+1},d^k)}{|\D(x^k,d^k)|}
=\frac{5/2-\zeta_k}{\zeta_k}\longrightarrow\infty.
\]
Thus sufficient descent fails, consistently with
Theorem~\ref{thm:dy}(b), and no fixed $c_3$ can satisfy the upper curvature bound
along this sequence.

\section{Conclusions}\label{sec:final}

The results of this paper clarify the distinct roles of descent
and line search conditions in the convergence of vector CG methods.
We have established global convergence of the pure FR and DY
methods under sufficient descent and the standard Wolfe conditions.
For FR, a suitable restriction on the curvature tolerance
guarantees sufficient descent. Under a generalized Wolfe line
search that includes the strong Wolfe conditions, descent alone
suffices for convergence. For DY, the same generalized Wolfe
conditions guarantee sufficient descent and global convergence.

The counterexamples show where the classical scalar properties
fail in vector optimization. The pure FR method may lose descent
under the strong Wolfe conditions, even for strongly convex
quadratic objectives and a curvature parameter smaller than $1/2$.
Moreover, both pure methods may fail to converge under the standard
Wolfe conditions, even when every search direction is a
$K$-descent direction, the initial level set is compact, and the
Jacobian is globally Lipschitz continuous.
Together, Theorems~\ref{thm:fr} and~\ref{thm:dy} and
Examples~\ref{ex:frcompact} and~\ref{ex:dy} answer (Q2) and (Q3).

Table~\ref{tab:summary} compares the scalar and vector results
under the assumptions of the cited results.
The conditions \eqref{eq:wmu} and \eqref{eq:gw} in the table
refer to restrictions on the upper tolerance in the Wolfe
conditions \eqref{eq:w1}--\eqref{eq:w2}.

\begin{table}[H]
\centering
\caption{Descent and global convergence of the pure
FR and DY methods.}
\label{tab:summary}
\small
\renewcommand{\arraystretch}{1.15}
\begin{tabularx}{\linewidth}{@{}
>{\raggedright\arraybackslash}p{0.25\linewidth}
>{\raggedright\arraybackslash}p{0.28\linewidth}
>{\raggedright\arraybackslash}X@{}}
\toprule
Method and conditions & Scalar case & Vector case\\
\midrule
FR, strong Wolfe,\newline $c_2<\tfrac12$
  & sufficient descent and convergence \cite{AlB85}
  & descent may fail (Example~\ref{ex:fr})\\
\addlinespace
FR, \eqref{eq:wmu}
  & sufficient descent and convergence \cite{DY96}
  & sufficient descent and convergence
    (Theorem~\ref{thm:fr}(a))\\
\addlinespace
FR, sufficient descent, standard Wolfe
  & convergence \cite{DY96}
  & convergence (Theorem~\ref{thm:fr}(b))\\
\addlinespace
FR, descent, standard Wolfe, $\mathcal L$ bounded
  & convergence \cite[Corollary~2]{Dai99}
  & convergence may fail (Example~\ref{ex:frcompact})\\
\addlinespace
FR, descent, \eqref{eq:gw}
  & convergence \cite{DY96}
  & convergence (Theorem~\ref{thm:fr}(c))\\
\addlinespace
DY, standard Wolfe
  & descent and convergence \cite{DY99}
  & descent (Theorem~\ref{thm:dy}),
    convergence may fail even with compact $\mathcal L$
    (Example~\ref{ex:dy})\\
\addlinespace
DY, \eqref{eq:gw}
  & sufficient descent and convergence \cite{DY99}
  & sufficient descent and convergence
    (Theorem~\ref{thm:dy}(a))\\
\bottomrule
\end{tabularx}
\end{table}

Our convergence analysis uses the finite total decrease of the
objective to control the growth of the search directions.
The resulting estimates provide a common framework for convergence
under the full FR and DY parameter bounds and also apply to
several scaled and hybrid methods.
The convergence results do not require a bounded level set.
Moreover, the same lower-bound assumption used in the convergence
analysis guarantees the existence of admissible Wolfe step sizes
for general ordering cones.

The estimate in Lemma~\ref{lem:key} does not rely on
a particular formula for the conjugate parameter.
It may therefore provide a useful tool for analyzing
other vector CG methods, especially those whose current
convergence results require strict fractions of the
conjugate parameters or additional assumptions on
the growth of the search directions.
Extending this estimate to the Riemannian setting
of \cite{CFS25,NH23} is another natural direction
for future research.

\section*{Funding}
L. F. Prudente was partially supported by CNPq,
Grants 307294/2023-4 and 403197/2025-2.

%

\section*{Use of generative AI}
AI tools (Claude, Anthropic, and Codex, OpenAI) were used
to check proofs and computations, to suggest alternative
arguments and simplifications, including the constructions
in Section~\ref{sec:nonconvergence}, and to revise the exposition.
All results were verified by the authors.
The authors assume responsibility for all content.

\begingroup
\small
\bibliographystyle{abbrv}
\bibliography{paper_FR_DY}
\endgroup

\end{document}